\documentclass[11pt,twoside,a4paper]{article}

\usepackage{graphicx,url}
\usepackage{setspace}
\usepackage{enumerate,lineno,setspace,float}
\usepackage[small,compact]{titlesec}
\usepackage{amsmath,amsfonts,amssymb,amsthm}
\usepackage{geometry}
\usepackage{fancyhdr}
\usepackage{amssymb,amsmath}
\usepackage{latexsym}
\usepackage{graphicx,float}
\usepackage{url}
\usepackage[capitalize]{cleveref} 
\usepackage{tikz}
\usepackage{caption}
\usepackage[numbers,sort&compress]{natbib}
\usepackage{mathrsfs}

\newtheorem{theorem}{Theorem}[section]
\newtheorem{lemma}[theorem]{Lemma}

\newtheorem{example}[theorem]{Example}

\newtheorem{observation}[theorem]{Observation}

\newtheorem{notation}[theorem]{Notation}
\usepackage{authblk}

\input amssym.def
\input amssym.tex

\newcommand{\be}{\begin{equation}}
\newcommand{\ee}{\end{equation}}
\newcommand{\ben}{\begin{equation*}}
\newcommand{\een}{\end{equation*}}
\newcommand{\bea}{\begin{eqnarray}}
\newcommand{\eea}{\end{eqnarray}}
\newcommand{\bean}{\begin{eqnarray*}}
\newcommand{\eean}{\end{eqnarray*}}
\def\non{\nonumber}

\def\qed{\hfill$\Box$\vspace{11pt}}

\def\diam{{\rm diam}}
\def\Span{{\rm span}}
\def\rn{{\rm rn}}

\def\C{\mathscr{C}}

\def\dis{\displaystyle}

\def\R{\mathcal{R}}
\def\({\left(}
\def\){\right)}
\def\[{\left[}
\def\]{\right]}

\begin{document}

\title{Optimal radio labelings for the Cartesian product of a path and the Heawood graph}

\author[a]{Payal Vasoya\thanks{E-mail: \texttt{prvasoya92@gmail.com}}}
\author[a]{Miral Doshi\thanks{E-mail: \texttt{doshi.miral55@gmail.com}}}
\author[b]{Devsi Bantva\thanks{E-mail: \texttt{devsi.bantva@gmail.com}}}
\affil[a]{{\small Gujarat Technological University, Ahmedabad 382 424, Gujarat, India}}
\affil[b]{{\small Lukhdhirji Engineering College, Morvi 363 642, Gujarat, India}}


\pagestyle{myheadings}
\markboth{\centerline{Payal Vasoya, Miral Doshi and Devsi Bantva }}{\centerline{Optimal radio labelings for the Cartesian product of a path and the Heawood graph}}

\date{}
\openup 0.8\jot
\maketitle

\begin{abstract} For a simple graph $G$, let $d(u,v)$ denote the distance between two vertices $u$ and $v$ in $G$, and let $\diam(G)$ denote the diameter of $G$. A \textit{radio labeling} of $G$ is a mapping  $f : V(G) \rightarrow \mathbb{N} \cup \{0\}$ such that $|f(u)-f(v)|\geq \diam(G) + 1 - d(u,v)$ holds for every pair of distinct vertices $u, v \in G$. The \textit{span} of $f$ is the difference between the largest and smallest integers assigned by $f$, and the \textit{radio number} $\rn(G)$ is the minimum span taking over all radio labelings of $G$. Let $P_n$ denote the path on $n$ vertices and let $H$ denote the Heawood graph (see Figure-1). In this paper we determine the radio number for the Cartesian product of a path and the Heawood graph, and characterize all optimal radio labelings of $P_n \Box H$. Moreover, we show that the radio number of $P_n \Box H$ exceed the lower bound given for $\rn(G)$ in \cite{Das}, where $G$ is any graph, and that this gap is arbitrarily large. \\


\end{abstract}

\section{Introduction}

The radio labelings originated from the classical \textit{channel assignment problem}, which concerns the allocation of limited frequency bands to a set of transmitters while avoiding signal interference. Two transmitters located close to one another must operate on sufficiently different frequencies to prevent overlapping signals, whereas transmitters that are farther apart can reuse frequencies that are closer in value. Modeling transmitters as vertices of a graph provides a natural mathematical abstraction of this setting: vertices represent transmitters, and edges connect pairs that are geographically close or prone to interference. A frequency assignment is then represented by a labeling of the vertices by nonnegative integers, and the constraint on channel separation can be expressed through distance requirements in the graph.  

In graph-theoretic terms, these considerations lead to the notion of a \textit{radio labeling} introduced by Chartrand \textit{et al.} \cite{Chartrand1,Chartrand2}. A \emph{radio labeling} of a graph $G = (V(G), E(G))$ is a function $f : V(G) \rightarrow \mathbb{N}_0$ ($\mathbb{N}_0$ is the set of positive integers including $0$) satisfying 
\begin{equation}\label{eqn:rn}
|f(u)-f(v)| \geq \diam(G)+1-d(u,v)
\end{equation}
for every pair of distinct vertices $u, v \in V(G)$, where $d(u,v)$ denote the distance between $u$ and $v$, and $\diam(G)$ denote the diameter of $G$. The assigned integer $f(u)$ is called the label of $u$ under $f$, and the span of $f$ is defined as $\Span(f) = \max\{|f(u)-f(v)| : u,v \in V(G)\}$. The radio number of $G$, denoted by $\rn(G)$, is defined as
$$\rn(G) := \min_{f}\{\Span(f)\}$$
where minimum taken over all radio labelings $f$ of $G$. A radio labeling $f$ that achieve this minimum is referred to as an \textit{optimal radio labeling}.

It is, without loss of generality, customary to assume some vertex is labeled $0$ in an optimal labeling; then the span equals the maximum label assigned. Observe that if $f$ is a radio labeling of a graph $G$, then $\Span(f) - f$ is also a radio labeling of $G$. It is clear from the definition that any radio labeling assigns different labels to distinct vertices of a graph $G$.

In this paper, we determine the radio number of the Cartesian product of a path and the Heawood graph, and we characterize all optimal radio labelings of the Cartesian product of a path and the Heawood graph. We show that the radio number of the Cartesian product of a path and the Heawood graph exceeds the lower bound given in \cite{Das}, and that the difference between the lower bound for the radio number of $P_n \Box H$ obtained using the technique given in \cite{Das} and the actual radio number of $P_n \Box H$ is arbitrarily large. 

\section{Preliminaries}

In this section, we provide the necessary terms and notations used for the present work. We follow \cite{West} for well-known graph-theoretic definitions and notations. For a simple connected graph $G = (V(G),E(G))$, the neighborhood $N(v)$ of a vertex $v \in V(G)$ is the set of all vertices adjacent to $v$. The distance $d(u,v)$ between two vertices $u$ and $v$ is the length of a shortest path joining $u$ and $v$. The diameter $\diam(G)$ of $G$ is $\max\{d(u,v) : u, v \in V(G)\}$. The eccentricity of $u \in V(G)$ is $\epsilon(u) = \max\{d(u,v) : v \in V(G)\}$. The center $C(G)$ of a graph $G$ is the subgraph induced by the vertices of minimum eccentricity. For a subset $S \subset V(G)$ and a vertex $u \in V(G) \setminus S$, define $d(u,S) = \min\{d(u,v) : v \in S\}$. A component of a graph $G$ is a maximal connected subgraph. A path $P_n$ is a graph with vertex set $V(P_n) = \{u_1,\ldots,u_n\}$ and edge set $E(P_n) = \{u_iu_{i+1} : 1 \leq i \leq n-1\}$. We fix the aforementioned vertex set of the path $P_n$ throughout this work. It is well known that $C(P_n)$ consists of a single vertex when $n$ is odd, and two adjacent vertices when $n$ is even. More precisely, with the above notations, $C(P_n) = \{u_{(n+1)/2}\}$ when $n$ is odd, and $C(P_n) = \{u_{n/2},u_{n/2+1}\}$ when $n$ is even. We view a path $P_n$ rooted at its center $C(P_n)$; if $C(P_n) = \{u_{(n+1)/2}\}$, then we view a path $P_n$ rooted at $u_{(n+1)/2}$, and if $C(P_n) = \{u_{n/2},u_{n/2+1}\}$, then we view a path $P_n$ rooted at both $u_{n/2}$ and $u_{n/2+1}$ (that is, both $u_{n/2}$ and $u_{n/2+1}$ are at level 0). Define a \textit{level function} $L : V(P_n) \rightarrow \{0,1,2,\ldots\}$ as follows:
$$L(u) = \min\{d(u,w) : w \in V(C(P_n))\} \mbox{ for each } u \in V(P_n),$$
and let
$$L(P_n) = \sum_{u \in V(P_n)}L(u).$$

A straightforward computation yields
\begin{equation*}
L(P_n) = \left\{
\begin{array}{ll}
\frac{1}{4}(n^2-1), & \mbox{if $n$ is odd}, \\[0.3cm]
\frac{1}{4}n(n-2), & \mbox{if $n$ is even}.
\end{array}
\right.
\end{equation*}



The Heawood graph $H$ is the graph with $V(H) = \{v_1,v_2,\ldots,v_{14}\}$ and $E(H)$ = \{$v_iv_{i+1}$, $v_1v_{14}$, $v_1v_6$, $v_2v_{11}$, $v_3v_{8}$, $v_4v_{13}$, $v_5v_{10}$, $v_7v_{12}$, $v_9v_{14}$ : $1 \leq i \leq 13$\} (see Figure-1). It is clear that $\diam(H) = 3$. We fix the above-mentioned vertex set for the Heawood graph throughout this work. Partition $V(H)$ into parity classes $A=\{v_i: i$ is odd$\}$ and $B=\{v_i: i$ is even$\}$; We refer to these as the \textit{categories} $H$. Two vertices $v_i$ and $v_j$ are in the same category if $i$ and $j$ have the same parity; otherwise, they are in the \emph{opposite} category.

\begin{observation}\label{Obs1} Let $H$ be the Heawood graph. Then the following statements hold:
\begin{enumerate}[\rm (a)]
\item If $u,v \in V(H)$ and $d(u,v) = 3$, then for any other vertex $w$, $d(u,w) \leq 2$ or $d(v,w) \leq 2$,
\item For three distinct vertices $u,v,w \in V(H)$, if $d(u,v) = d(v,w) = 3$, then $d(u,w) = 2$, 
\item Let $u \in V(H)$. Then, among the seven vertices in the opposite category of $u$, exactly three are adjacent to $u$, and the remaining four are at distance three from $u$.
\end{enumerate}
\end{observation}

Let $G_1=(V(G_1),E(G_1))$ and $G_2=(V(G_2),E(G_2))$ be two graphs. The Cartesian product of $G_1$ and $G_2$ is the graph $G_1 \Box G_2$, with $V(G_1 \Box G_2) = V(G_1) \times V(G_2)$,  such that two vertices $(a_1,b_1)$ and $(a_2,b_2)$ are adjacent if $a_1=a_2$ and $(b_1,b_2) \in E(G_2)$, or $b_1=b_2$ and $(a_1,a_2) \in E(G_1)$. It is clear from the definition that $d_{G_1 \Box G_2} ((a_1,b_1), (a_2,b_2)) = d_{G_1} (a_1,a_2) + d_{G_2} (b_1,b_2)$.

\begin{observation} For $P_n \Box H$,
$|V(P_n \Box H)| = 14n$ and $\diam(P_n \Box H) = n+2$.
\end{observation}

Let $V(P_n \Box H) = (x_0,x_1,\ldots,x_{14n-1})$. Note that each $x_t\;(0 \leq t \leq 14n-1)$ is an ordered pair $(u_{i_t},v_{j_t})$, with $u_{i_t} \in V(P_n)$ and $v_{j_t} \in V(H)$, where $1 \leq i_t \leq n, 1 \leq j_t \leq 14$. Define a level function $L : V(P_n \Box H) \rightarrow \mathbb{N} \cup \{0\}$ by
$$L(x_t) := L(u_{i_t}) = \min\{d(u_{i_t},w) : w \in V(C(P_n))\}\mbox{ for each } x_t = (u_{i_t},v_{j_t}) \in V(P_n \Box H),$$
and define
$$L(P_n \Box H)= \sum_{x_t \in V(P_n \Box H)}L(x_t).$$

Observe that
\begin{equation*}
L(P_n \Box H) = 14L(P_n) = \left\{
\begin{array}{ll}
\frac{7}{2}(n^2-1), & \mbox{if $n$ is odd}, \\[0.3cm]
\frac{7}{2}n(n-2), & \mbox{if $n$ is even}.
\end{array}
\right.
\end{equation*}

Define the \emph{layers} of $P_n \Box H$ by $S_k,k=1,2,\ldots,n$, where $S_k = \{x_t : x_t = (u_{i_t},v_{j_t}) \mbox{ such that } i_t =k\}$. The subgraph induced by each layer $S_k$ is isomorphic to the Heawood graph $H$. We partition $V(P_n \Box H)$ into disjoint sets by putting  vertices on the same side in the same set. If $n$ is even, there are two sides, namely, 
$$X_1 = \displaystyle\bigcup_{i=1}^{\frac{n}{2}} S_i\;\; \mbox{ and }\;\;  X_2=\displaystyle\bigcup_{i=\frac{n}{2}+1}^{n}S_i.$$ 
If $n$ is odd, there are three sides, namely, $$X_1=\displaystyle\bigcup_{i=1}^{\frac{n-1}{2}}S_i,\;\; X_2=S_{\frac{n+1}{2}}\;\; \mbox{ and }\;\; X_3=\bigcup_{i=\frac{n+3}{2}}^{n}S_i.$$ Similarly, we partition $V(P_n \Box H)$ into disjoint sets by putting vertices on the same category in the same set. For any $n$, $P_n \Box H$ has two categories, namely, 
$$Y_1 = \{(u_{i_t},v_{j_t}) : \mbox{$j_t$ is odd}\} \mbox{ and }Y_2 = \{(u_{i_t},v_{j_t}) : \mbox{$j_t$ is even}\}.$$ 
It is clear that $|Y_1| = |Y_2| = 7n$. Two vertices, $x_a = (u_{i_a},v_{j_a})$ and $x_b = (u_{i_b},v_{j_b})$ of $V(P_n \Box H)$, are called
on the \emph{opposite} or on the \emph{same} side depending on whether $x_a$ and $x_b$ are on opposite or the same side of $P_n \Box H$, and in \emph{opposite} or \emph{same} category depending on whether $v_{j_a}$ and $v_{j_b}$ are in opposite or in the same category in $H$. 

For any $x_a, x_b \in V(P_n \Box H)$, where $x_a = (u_{i_a},v_{j_a})$ and $x_b = (u_{i_b},v_{j_b})$, define
\begin{equation*}
\phi(x_a, x_b) := \left\{
\begin{array}{ll}
\min\{L(x_a),L(x_b)\}, & \mbox{ if $x_a$ and $x_b$ are on the same side}, \\
0, & \mbox{ otherwise},
\end{array}
\right.
\end{equation*}
\begin{equation*}
\delta(x_a,x_b) := \left\{
\begin{array}{ll}
1, & \mbox{ if $n$ is even, and $x_a$ and $x_b$ are on the opposite sides}, \\
0, & \mbox{ otherwise},
\end{array}
\right.
\end{equation*}
\begin{equation*}
d_1(x_a,x_b) := d(u_{i_a},u_{i_b}) \mbox{ in $P_n$ and } d_2(x_a,x_b) := d(v_{j_a},v_{j_b})\mbox{ in $H$}.\hspace{1.5cm}
\end{equation*}

\begin{lemma}\label{lem:phi} For $x_a,x_b \in V(P_n \Box H)$, the following hold.
\begin{enumerate}[\rm (a)]
\item $0 \leq \phi(x_a,x_b) \leq \lfloor (n-1)/2 \rfloor$;  \item $\phi(x_a,x_b) = 0$ if and only if (i) $x_a$ and $x_b$ are on opposite sides of $P_n \Box H$; (ii) at least one of $x_a$ and $x_b$ is in $S_{(n+1)/2}$ when $n$ is odd, or in one of $S_{n/2}$ and $S_{n/2+1}$ when $n$ is even; 
\item if $n$ is odd and $x_a,x_b \in V(P_n \Box H) \setminus S_{(n+1)/2}$ are on the same side, then $\phi(x_a,x_b) \geq 1$; 
\item $0 \leq d_1(x_a,x_b) \leq n-1$ and $0 \leq d_2(x_a,x_b) \leq 3$;
\item the distance between $x_a$ and $x_b$ is given by
\begin{equation}\label{eqn:dxaxb}
d(x_a,x_b) = L(x_a)+L(x_b)+\delta(x_a,x_b)-2\phi(x_a,x_b)+d_2(x_a,x_b).
\end{equation}
\end{enumerate}
\end{lemma}

\begin{notation}Let $\vec{V}= (x_0,x_1,\ldots,x_{14n-1})$ be an ordering of $V(P_n \Box H)$. For a triple $(x_a,x_b,x_c)$ consecutive vertices set, denote $$\ell(x_a,x_b,x_c) = \min\{d_1(x_a,x_b), d_1(x_b,x_c)\}.$$
\end{notation}

\noindent A triple $(x_l,x_{l+1},x_{l+2})$ of $\vec{V}$ is called:
\begin{itemize}
\item \emph{Proper} if $d_2(x_t,x_{t+1}) = 3$ for $t=l,l+1$, and $d_2(x_l,x_{l+2}) = 2$;
\item \textit{Semi-proper} if $\min\{d_2(x_l,x_{l+1})$, $d_2(x_{l+1},x_{l+2})\} = 2$, $\max\{d_2(x_l,x_{l+1}),d_2(x_{l+1},x_{l+2})\}=3$ and $d_2(x_l,x_{l+2})=3$;
\item \textit{Weak-proper} if $d_2(x_t,x_{t+1}) = 2$ for $t=l,l+1$ and $d_2(x_l,x_{l+2})=2$.
\end{itemize}    
A segment of vertices $(x_k,\ldots,x_m)$ (possibly with $k=0$, $m=14n-1$) in the ordering $\vec{V} = (x_0,x_1,\ldots,x_{14n-1})$ of $V(P_n \Box H)$ is called \textit{proper}, \textit{semi-proper} and \textit{weak-proper} if each subsegment $(x_{k+t},x_{k+t+1},x_{k+t+2})$ of vertices is proper, semi-proper and weak-proper for $t=0,2,\ldots,m-k-2$. A triple $(x_l,x_{l+1},x_{l+2})$ in an ordering  $\vec{V}_f : (x_0,x_1,\ldots,x_{14n-1})$ is called  \textit{smooth} if $\{x_l,x_{l+1},x_{l+2}\} \cap S_x \neq \emptyset$ for each $x \in \{1,(n+1)/2,n\}$ when $n$ is odd. An ordering $\vec{V}_f = (x_0,x_1,\ldots,x_{14n-1})$ of $V(P_n \Box H)$, where $n$ is odd, is called a \textit{smooth ordering} if the vertices of $S_x$ for $x = 1,(n+1)/2,n$ are ordered exclusively as smooth triples. 

\begin{observation}\label{Obs2} Let $\vec{V} = (x_0,x_1,\ldots,x_{14n-1})$ be an ordering of $V(P_n \Box H)$ and  $(x_k,\ldots,x_m)$ be a segment with (i) $x_t$ and $x_{t+1}$ are on opposite side for $k \leq t \leq m-1$, and (ii) $x_t$ and $x_{t+2}$ are on the same side for $k \leq t \leq m-2$. Then the following statements hold:
\begin{enumerate}[\rm (a)]
\item If $(x_k,\ldots,x_m)$ is proper, then $x_l$ and $x_{l+1}$ are of opposite category, and $x_l$ and $x_{l+2}$ of same category for $k \leq l \leq m-1$ (that is, the vertices ordered from $(x_k,\ldots,x_{m})$ on same side are of same category and on opposite side are of opposite category).
\item If $(x_k,\ldots,x_m)$ is semi-proper, then $x_l$ and $x_{l+2}$ are of opposite category for $k \leq l \leq m-1$ (that is, the vertices ordered from $(x_k,\ldots,x_m)$ on same side are of opposite category).
\item If $(x_k,\ldots,x_m)$ is weak-proper, then $x_l, x_{l+1}, x_{l+2}$ all are of the same category for $k \leq l \leq m-1$ (that is, all the vertices ordered from $(x_k,\ldots,x_{m})$ are of the  same category).
\end{enumerate}
\end{observation}

\begin{lemma}\label{lem:max} Let $\vec{V} = (x_0,x_1,\ldots,x_{14n-1})$ be an ordering of $V(P_n \Box H)$, where $n$ is even, such that $x_t$ and $x_{t+1}$ are on opposite sides for all $0 \leq t \leq 14n-2$. Then there exists at least one $t$ with $d_2(x_t,x_{t+1})=2$. Moreover, if there is exactly one such pair then its index satisfies $t=7n-1$.
\end{lemma}
\begin{proof} Suppose to the contrary that $d_2(x_t,x_{t+1})=3$ for all $0 \leq t \leq 14n-2$. Since $n$ is even, recall that $P_n \Box H$ has two sides, $X_1$ and $X_2$, and two categories, $Y_1$ and $Y_2$. Moreover, note that $|X_1 \cap Y_1| = |X_1 \cap Y_2| = |X_2 \cap Y_1| = |X_2 \cap Y_2|= 7n/2$. By Observation \ref{Obs2}, if consecutive vertices $x_t$ and $x_{t+1}$ always have $d_2(x_t,x_{t+1}) = 3$, then the vertices appearing on a given side in the ordering $\vec{V}$ would all belong to the same category, contradicting the balanced distribution of categories across sides. Hence there must be at least one index $t$ with $t$ with $d_2(x_t,x_{t+1})=2$. 

If there is exactly one such index, by symmetry and counting of the four-equal sized subsets $X_1 \cap Y_1, X_1 \cap Y_2, X_2 \cap Y_1, X_2 \cap Y_2$ force the unique transition that switches category counts to occur at the midpoint of the ordering; that midpoint index is $7n-1$, which completes the proof. \qed
\end{proof}

We now introduce notation for radio labelings organized by the induced ordering. Let $f$ be any radio labeling of $P_n \Box H$, and let $\vec{V}_f:=(x_0,x_1,\ldots,x_{14n-1})$ be the ordering of $V(P_n \Box H)$ induced by $f$ such that $0 = f(x_0) < f(x_1) < \ldots < f(x_{14n-1}) = \Span(f)$. Denote by $\R$ the set of all radio labelings of $P_n \Box H$; denote by $\R_o$ the set of all radio labelings whose induced ordering satisfies that consecutive vertices lie on opposite sides and by $\R_s = \R \setminus \R_o$.

\section{Main results}
In this section, we determine the radio number of the Cartesian product of a path and the Heawood graph, and characterize all optimal radio labelings for it. The main result is as follows:
\begin{theorem}\label{thm:PHL} Let $n \geq 3$ be any integer. Then
\begin{equation}\label{rn:PHL}
\rn(P_{n} \Box H) = \left\{
\begin{array}{ll}
68, & \mbox{ if $n = 3$}, \\ [0.2cm]
7n^2-n+9, & \mbox{ if $n > 3$ is odd}, \\ [0.2cm]
n(7n+6), & \mbox{ if $n > 3$ is even},
\end{array}
\right.
\end{equation}
Moreover, an injective function $f : V(P_n \Box H) \rightarrow \mathbb{N} \cup \{0\}$ is an optimal radio labeling of $P_n \Box H$ if and only if the following all hold:
\begin{enumerate}[\rm (a)]
\item The ordering $\{x_0,x_1,\ldots,x_{14n-1}\}$ of $V(P_n \Box H)$ induced by $f$ such that $0 = f(x_0) < f(x_1) < \cdots < f(x_{1n-1}) = \Span(f)$ satisfies the following conditions: \\
    \textsf{Case-1:} $n$ is odd.
\begin{enumerate}[\rm (1)]
\item $x_t$ and $x_{t+1}$ are on opposite sides for every $0 \leq t \leq 14n-1$,
\item one of the following holds:
\begin{enumerate}[\rm (i)]
\item $x_0 \in S_{(n+1)/2}$ and $x_{14n-1} \in N(x_0), x_{14n-1} \notin S_{(n+1)/2}$ and $(x_t,x_{t+1},x_{t+2})$ is proper for all $t=0,2,\ldots,14n-4$ with $\ell(x_t,x_{t+1},x_{t+2}) = \lfloor n/2 \rfloor$ except one triple which is semi-proper with $\ell(x_t,x_{t+1},x_{t+2}) = \lfloor n/2 \rfloor+1$;
\item $x_0,x_{14n-1} \in S_{(n+1)/2}$ and two triples among $(x_t,x_{t+1},x_{t+2})$ for $t=0,2,\ldots,14n-4$ are semi-proper and all others are proper, and $\ell(x_t,x_{t+1},x_{t+2}) = \lfloor n/2 \rfloor$,
        \end{enumerate}
        \item $d_2(x_{14n-2},x_{14n-1})=3$.
    \end{enumerate}
    \textsf{Case-2:} $n$ is even.
    \begin{enumerate}[\rm (1)]
\item $x_0 \in S_{n/2}$ and $x_{14n-1} \in S_{n/2+1}$;
\item $x_t$ and $x_{t+1}$ are on opposite side for every $0 \leq t \leq 14n-2$;
\item for $t = 0,2,\ldots,14n-4$, $(x_t,x_{t+1},x_{t+2})$ is proper or semi-proper with $\ell(x_t,x_{t+1},x_{t+2}) = n/2$ or there are equal number of proper triples $(x_t,x_{t+1},x_{t+2})$ with $\ell(x_t,x_{t+1},x_{t+2}) = n/2-1$ and semi-proper triples $(x_t,x_{t+1},x_{t+2})$ with $\ell(x_t,x_{t+1},x_{t+2}) = n/2+1$;
        \item $d_2(x_{14n-2},x_{14n-1})=3$.
    \end{enumerate}
  \item $f$ is defined by $f(x_0) = 0$ and $f(x_{t+1}) = f(x_t)+d+1+J_t-d(x_t,x_{t+1})$ for $0 \leq t \leq 14n-2$, where $J_t$ are the integers defined as follows:
\begin{equation}\label{eqn:Jt}
J_t = \left\{
\begin{array}{ll}
2, & \mbox{ if $n$ is odd and $(x_{t-1},x_{t},x_{t+1})$ is proper with $l(x_{t-1},x_{t},x_{t+1}) = (n+1)/2$}, \\
1, & \mbox{ if $n$ is even and $(x_{t-1},x_{t},x_{t+1})$ is proper with $l(x_{t-1},x_{t},x_{t+1}) = n/2$}, \\
0, & \mbox{ otherwise}.
\end{array}
\right.
\end{equation}
\end{enumerate}
\end{theorem}

Since our main result is related to the Cartesian product of a path and the Heawood graph, we first present the radio numbers of these two underlying graph families for completeness. The radio number of paths is determined by Liu and Zhu in \cite{Liu}.

\begin{theorem}(\cite{Liu}) Let $n \geq 4$ be an integer. Then
\begin{equation*}
\rn(P_n) = \left\{
\begin{array}{ll}
2m^2+2, & \mbox{ if } n=2m+1, \\
2m(m-1)+1, & \mbox{ if }n=2m. \\
\end{array}
\right.
\end{equation*}
\end{theorem}

\begin{theorem}\label{thm:H} For the Heawood graph $H$, $\rn(H) = 13$.
\end{theorem}
\begin{proof} Note that any radio labeling $f$ is a one-to-one integral function from the vertex set $V(G)$ to the set of non-negative integers, and hence $\Span(f) \geq n-1$ for any graph $G$ on $n$ vertices. Since $H$ has $14$ vertices, $\rn(H) \geq 13$. We now prove that this bound is tight by giving a radio labeling $f$ whose span is equal to $13$. Denote each vertex $v_i \in V(H)$ by $v_{3j+t}$ where $i=3j+t$ for some $j \in \{0,1,2,\ldots\}$ and $t \in \{1,2,3\}$. Define a labeling by $f(v_{3j+t}) = j+5t-5$. See Figure \ref{Fig1}. One can easily verify that $f$ is a radio labeling of $H$ whose span is equal to $13$. \qed
\begin{figure}[ht!]
    \centering
    \includegraphics[width=2.5in]{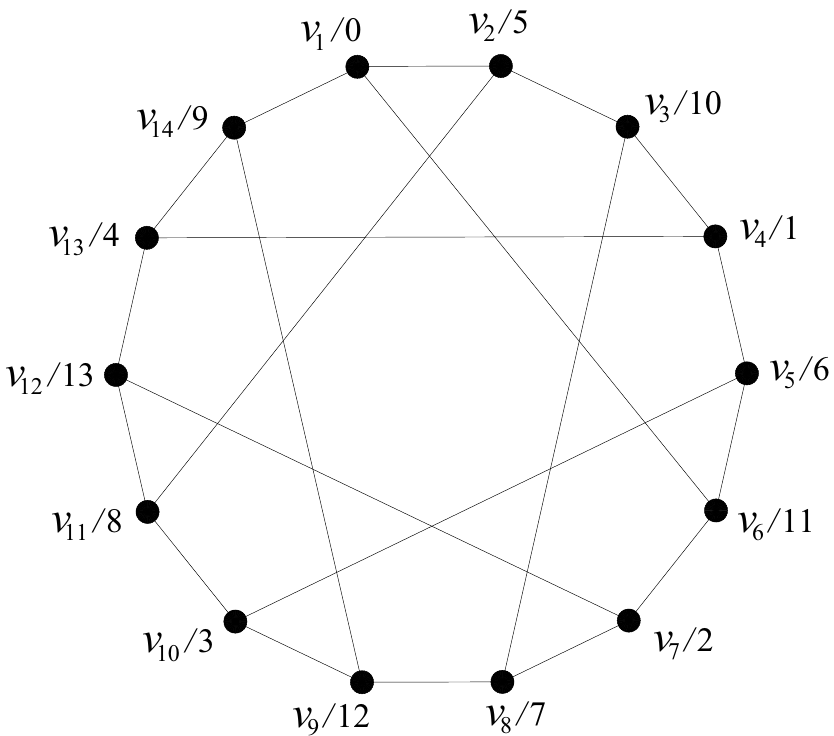}
    \caption{\rn(H) = 13.}
    \label{Fig1}
\end{figure}
\end{proof}

Next, we proceed to prove our main result. We need some preparation in order to determine the radio number of Cartesian product of a path and the Heawood graph.

\begin{lemma}\label{Prop:2} Let $G$ be a simple graph of order $p$ and diameter $d \geq 2$. Let $\vec{V} = (x_0, x_1, \ldots, x_{p-1})$ be an ordering of $V(G)$. For $0 \leq t \leq p-2$, denote $d_t = d-d(x_t,x_{t+1}) \geq 0$ and define a function $f$ by $f(x_0) = 0$ and $f(x_{t+1}) = 1+d_t+f(x_t)$ for $0 \leq t \leq p-2$. If for every pair of vertices $x_i$ and $x_j$ such that $j-i+\sum_{t=i}^{j-1}d_t < d$, the inequality \eqref{eqn:rn} holds, then  $f$ is a radio labeling of $G$.
\end{lemma}
\begin{proof} Let $x_i$ and $x_j$ be two vertices with $0 \leq i < j \leq p-1$. If $j-i+\sum_{t=i}^{j-1}d_t < d$, then the inequality \eqref{eqn:rn} is already satisfied by hypothesis. Assume $j-i+\sum_{t=i}^{j-1}d_t \geq d$. Then, 
\bean
f(x_j)-f(x_i) & = & \sum_{t=i}^{j-1}(f(x_{t+1})-f(x_t)) = \sum_{t=i}^{j-1}(d+1-d(x_t,x_{t+1})) \\
& = & \sum_{t=i}^{j-1}(1+d_t) = (j-i)+\sum_{t=i}^{j-1}d_t \geq d \geq d+1-d(x_i,x_j), \\
\eean
So, the inequality \eqref{eqn:rn} also holds in this case. Hence, $f$ is a radio labeling of $G$ as the inequality \eqref{eqn:rn} holds for all $x_i$ and $x_j,\; 0 \leq i < j \leq p-1$. \qed
\end{proof}

\begin{lemma}\label{lemdist} Let $n > 3$ be an even integer, and let $f \in \R_o$ be a radio labeling of $P_n \Box H$ whose induced ordering is $\vec{V}_f = (x_0,x_1,\ldots,x_{14n-1})$. Suppose a triple $(x_l,x_{l+1},x_{l+2})$ satisfies $\ell(x_l,x_{l+1},x_{l+2}) \leq n/2-\alpha$ for some integer $\alpha \geq 1$. Then, there exists another triple $(x_{l'},x_{l'+1},x_{l'+2})$ in the ordering with  $\ell(x_{l'},x_{l'+1},x_{l'+2}) \geq n/2+\alpha$.
\end{lemma}
\begin{proof} Assume to the contrary that every triple $(x_t,x_{t+1},x_{t+2})$ satisfies $\ell(x_t,x_{t+1},x_{t+2}) \leq n/2 + \alpha - 1$. Recall that $P_n \Box H$ has two sides, $X_1$ and $X_2$ with $|X_1| = |X_2|$, where $X_1 = \cup_{k=1}^{n/2}S_k$, $X_2 = \cup_{k=n/2+1}^{n}S_k$. Note that for any $1 \leq k \leq n/2$, $|S_k| = |S_{n/2+k}| = 14$, and $d_1(x_a,x_b) = n/2$ for $x_a \in S_k$ and $x_b \in S_{n/2+k}$. Moreover, for any $x_a = (u_{i_a},v_{j_a}), x_b = (u_{i_b},v_{j_b}) \in V(P_n \Box H)$, $d_1(x_a,x_b) = |i_a-i_b|$, where $i_a, i_b \in \{1,2,\ldots,n\}$. Now $d(x_t,x_{t+1}) \leq n/2-\alpha$ for $t=l,l+1$. Without loss of generality, assume that $x_l,x_{l+2} \in X_1$ and $x_{l+1} \in X_2$. This implies that $i_l,i_{l+2} \in \{1,\ldots,n/2\}$, and $i_{l+1} \in \{n/2,\ldots,n\}$. Again, without loss of generality, assume that $i_l \geq i_{l+2}$, then note that there exist a vertex $x_{l'+1} \in \cup_{t=i_l+n/2}^{n}S_t$ such that $d(x_{l'},x_{l'+1}) \geq n/2+\alpha$ and $d(x_{l'+1},x_{l'+2}) \geq n/2+\alpha$. That is, $\ell(x_{l'},x_{l'+1},x_{l'+2}) \geq n/2+\alpha$ which is a contradiction with our assumption. Thus, the proof is complete. \qed
\end{proof}

Let $f$ be any radio labeling of $P_n \Box H$ whose induced ordering is $\vec{V}_f=(x_0,x_1,\ldots,x_{14n-1})$.
Define
\begin{equation}\label{eqn:j1}
J_f(x_{i},x_{i+1}) = f(x_{i+1})-f(x_{i})+d(x_i,x_{i+1})-(d+1),\;\; 0 \leq i \leq 14n-2.
\end{equation}
It is clear that $J_f(x_i,x_{i+1}) \geq 0$ for all $0 \leq i \leq 14n-2$. Moreover, for $0 \leq i \leq n-2$,  $J_f(x_i,x_{i+1})=0$ if and only if $f$ is defined by $f(x_0)=0$ and $f(x_{i+1})=f(x_i)+d+1-d(x_i,x_{i+1})$. For any segment $(x_i,x_{i+1},\ldots,x_j)$ of the ordering
$\vec{V}_f=(x_0,x_1,\ldots,x_{14n-1})$ of $V(P_n \Box H)$, where $0 \leq i < j \leq 14n-1$, define
\begin{equation}\label{eqn:Jij}
J_f(x_i,\ldots,x_j) = \dis\sum_{t=i}^{j-1}J_f(x_t,x_{t+1})
\end{equation}
and
\begin{equation*}
J_f(P_n \Box H) := J_f(x_0,\ldots,x_{14n-1}) = \dis\sum_{i=0}^{14n-2}J_f(x_i,x_{i+1}).
\end{equation*}

\begin{lemma}\label{lem3} Let $\alpha \geq 0$ and $n \geq 3$ be two integers. Let $f \in \R_o$ be any radio labeling of $P_n \Box H$, and let $\vec{V}_f=(x_0,x_1,\ldots,x_{14n-1})$ be the ordering of $V(P_n \Box H)$ induced by $f$. Suppose that one of the following conditions holds for the triple  $(x_l,x_{l+1},x_{l+2})$ for some $l$, where $0 \leq l \leq 14n-3$:
\begin{enumerate}[\rm (a)]
\item $(x_l,x_{l+1},x_{l+2})$ is proper and $d_1(x_{t},x_{t+1}) \geq \lfloor n/2 \rfloor + \alpha$ for $t = l, l+1$;
\item $(x_l,x_{l+1},x_{l+2})$ is semi-proper or weak-proper and $d_1(x_t,x_{t+1}) \geq \lfloor n/2 \rfloor + 1 + \alpha$ for $t = l,l+1$.
\end{enumerate}
Then, 
\begin{equation}\label{eqn:jfx3}
J_f(x_l,x_{l+1},x_{l+2}) \geq \left\{
\begin{array}{ll}
2\alpha, & \mbox{ if $n$ is odd}, \\ [0.2cm]
2\alpha+1, & \mbox{ if $n$ is even}.
\end{array}
\right.
\end{equation}
\end{lemma}
\begin{proof} By \eqref{eqn:j1} and \eqref{eqn:Jij}, we have  
\begin{eqnarray}
J_f(x_l,x_{l+1},x_{l+2}) & = & J_f(x_l,x_{l+1})+J_f(x_{l+1},x_{l+2}) \nonumber \\
     & = & f(x_{l+2})-f(x_l)+d(x_l,x_{l+1})+d(x_{l+1}, x_{l+2})-2(d+1).  \label{eqn:j4}
\end{eqnarray}

\noindent
Since $f$ is a radio labeling, we have $f(x_{l+2})-f(x_l) \geq d+1-d(x_l,x_{l+2})$. Substituting this into  \eqref{eqn:j4}, we obtain 
\begin{eqnarray}
J_f(x_l,x_{l+1},x_{l+2}) & \geq & d(x_l,x_{l+1})+d(x_{l+1},x_{l+2})-d(x_l,x_{l+2})-d-1. \label{eqn:j5}
\end{eqnarray}

Since $f \in \R_o$, $x_t$ and $x_{t+1}$ are on the opposite sides for $t=l,l+1$. We proceed the proof by considering the following two cases. Recall that $d=n+2$.

\textsf{Case-1:} $n$ is odd. Then $d$ is odd as $d = n+2$. Assume that $x_l$ and $x_{l+2}$ are in the opposite sides, then one of the vertex in $\{x_l,x_{l+1},x_{l+2}\}$ is in $S_{(n+1)/2}$.  Therefore, we have $\min\{d_1(x_l,x_{l+1}), d_1(x_{l+1},x_{l+2})\} = (d-3)/2 = (n-1)/2$ and hence, by (a), we obtain $\alpha = 0$. Thus, $J_f(x_l,x_{l+1},x_{l+2}) \geq 0$.

Assume that $x_l$ and $x_{l+2}$ are on the same side. Then, by \eqref{eqn:dxaxb} and Lemma  \ref{lem:phi}, $d(x_t,x_{t+1}) =  L(x_t)+L(x_{t+1})+d_2(x_t,x_{t+1})$ for $t=l, l+1$, and $d(x_l,x_{l+2}) = L(x_l)+L(x_{l+2})-2\phi(x_l,x_{l+2})+d_2(x_l,x_{l+2})$. Substituting all these into \eqref{eqn:j5}, we have
\begin{equation}
J_f(x_{l},x_{l+1},x_{l+2}) \geq 2L(x_{l+1})+2\phi(x_l,x_{l+2})+d_2(x_l,x_{l+1})+d_2(x_{l+1},x_{l+2})-d_2(x_l,x_{l+2})-d-1.\label{eqn:j6}
\end{equation}

Suppose that (a) holds. Then, by Observation \ref{Obs1} and Lemma \ref{lem:phi}, $d_2(x_l,x_{l+1})$ = $d_2(x_{l+1},x_{l+2})$ = $3$ and
$d_2(x_l,x_{l+2}) \leq 2$, and $d_1(x_t, x_{t+1}) = L(x_t) + L(x_{t+1}) \geq (d-3)/2+\alpha$ for $t=l, l+1$. Since  $\phi(x_l,x_{l+2}) = \min\{L(x_l), L(x_{l+2})\}$, we have $L(x_{l+1})+\phi(x_l,x_{l+2}) \geq (d-3)/2+\alpha$. Substituting all these into \eqref{eqn:j6} we obtain, $J_f(x_l, x_{l+1}, x_{l+2}) \geq  2\alpha$.

Suppose that (b) holds. Assume that $(x_l,x_{l+1},x_{l+2})$ is semi-proper. Then, by Observation \ref{Obs1} and Lemma \ref{lem:phi}, $d_2(x_l,x_{l+1})+d_2(x_{l+1},x_{l+2})=5$ and $d_2(x_l,x_{l+2}) \leq 3$, and $d_1(x_t,x_{t+1}) \geq (d-1)/2+\alpha$ for $t=l,l+1$. Since $\phi(x_l,x_{l+2})=\min\{L(x_l),L(x_{l+2})\}$, we have $L(x_{l+1})+\phi(x_l,x_{l+2}) \geq (d-1)/2+\alpha$. Substituting all these into \eqref{eqn:j6} we obtain, $J_f(x_l,x_{l+1},x_{l+2}) \geq 2\alpha$.

Assume that $(x_l,x_{l+1},x_{l+2})$ is weak-proper. Then, by Observation \ref{Obs1} and Lemma \ref{lem:phi}, $d_2(x_l,x_{l+1}) = d_2(x_{l+1},x_{l+2})=2$ and $d_2(x_l,x_{l+2}) \leq 2$, and $d_1(x_t,x_{t+1}) \geq (d-1)/2+\alpha$ for $t=l,l+1$. Since $\phi(x_l,x_{l+2})=\min\{L(x_l),L(x_{l+2})\}$, we have $L(x_{l+1})+\phi(x_l,x_{l+2}) \geq (d-1)/2+\alpha$. Substituting all these into \eqref{eqn:j6} we obtain, $J_f(x_l,x_{l+1},x_{l+2}) \geq 2\alpha$.

\textsf{Case-2:} $n$ is even. Then $d$ is even as $d=n+2$. By the assumption $f \in \R_o$, $x_l$ and $x_{l+2}$ are in the same side. By \eqref{eqn:dxaxb} and Lemma \ref{lem:phi}, $d(x_t,x_{t+1}) = L(x_t)+L(x_{t+1})+1+d_2(x_t,x_{t+1})$ for $t=l,l+1$ and $d(x_l,x_{l+2}) = L(x_l)+L(x_{l+2})-2\phi(x_l,x_{l+2})+d_2(x_l,x_{l+2})$. Substituting all these into \eqref{eqn:j5}, we have
\begin{equation}
J_f(x_l,x_{l+1},x_{l+2}) \geq 2L(x_{l+1})+2\phi(x_l,x_{l+2})+d_2(x_l,x_{l+1})+d_2(x_{l+1},x_{l+2})-d_2(x_l,x_{l+2})-d+1.\label{eqn:j7}
\end{equation}

Suppose that (a) holds. Then, by Observation \ref{Obs1} and Lemma \ref{lem:phi}, we have $d_2(x_l,x_{l+1}) = d_2(x_{l+1},x_{l+2}) = 3$ and
$d_2(x_l,x_{l+2}) \leq 2$, and  $d_1(x_t,x_{t+1})=L(x_t)+L(x_{t+1})+1 \geq (d-2)/2+\alpha$ for $t=l,l+1$. Since $\phi(x_l,x_{l+2}) = \min\{L(x_l),L(x_{l+2})\}$, we obtain $L(x_{l+1})+\phi(x_l,x_{l+2}) \geq (d-4)/2+\alpha$. Substituting all these in \eqref{eqn:j7} we obtain,
$J_f(x_l,x_{l+1},x_{l+2}) \geq  2\alpha+1$.

Suppose that (b) holds. Assume that $(x_l,x_{l+1},x_{l+2})$ is semi-proper. Then, by Observation \ref{Obs1} and Lemma \ref{lem:phi}, we have $d_2(x_l,x_{l+1})+d_2(x_{l+1},x_{l+2})=5$ and $d_2(x_l,x_{l+2}) \leq 3$, and $d_1(x_t,x_{t+1}) = L(x_t)+L(x_{t+1})+1 \geq d/2+\alpha$ for $t=l,l+1$. Since $\phi(x_l,x_{l+2}) = \min\{L(x_l),L(x_{l+2})\}$, we obtain $L(x_{l+1})+\phi(x_l,x_{l+2}) \geq d/2-1+\alpha$. Substituting all these in \eqref{eqn:j7} we obtain, $J_f(x_l,x_{l+1},x_{l+2}) \geq 2\alpha+1$.

Assume that $(x_l,x_{l+1},x_{l+2})$ is weak-proper. Then, by Observation \ref{Obs1} and Lemma \ref{lem:phi}, we have $d_2(x_l,x_{l+1}) = d_2(x_{l+1},x_{l+2}) = 2$ and $d_2(x_l,x_{l+2}) \leq 2$, and $d_1(x_t,x_{t+1}) = L(x_t)+L(x_{t+1})+1 \geq d/2+\alpha$ for $t=l,l+1$. Since  $\phi(x_l,x_{l+2}) = \min\{L(x_l),L(x_{l+2})\}$, we obtain $L(x_{l+1})+\phi(x_l,x_{l+2}) \geq d/2-1+\alpha$. Substituting all these in \eqref{eqn:j7} we obtain, $J_f(x_l,x_{l+1},x_{l+2}) \geq 2\alpha+1$. \qed
\end{proof}

\begin{lemma}\label{lemjG} Let $f \in \R_o$ be any radio labeling of $P_n \Box H$, and let $\vec{V}_f=(x_0,x_1,\ldots,x_{14n-1})$ be the ordering of $V(P_n \Box H)$ induced by $f$, which is proper except one triple which is semi-proper when $n$ is even. Then
\begin{equation}\label{jump:PnH}
J_f(P_n \Box H) \geq
\left\{
\begin{array}{ll}
0, & \mbox{ if $n$ is odd}, \\
7n-2, & \mbox{ if $n$ is even}.
\end{array}
\right.
\end{equation}
\end{lemma}
\begin{proof} Assume that $n$ is odd. Then, observe that the result follows, as $J_f(x_t,x_{t+1}) \geq 0$ for all $0 \leq t \leq 14n-2$.

Assume that $n$ is even. In this case, recall that $P_n \Box H$ has two sides, $X_1$ and $X_2$ with $|X_1| = |X_2|$, where $X_1 = \cup_{k=1}^{n/2}S_k, X_2 = \cup_{k=(n/2)+1}^{n}S_k$. Moreover, for any $1 \leq k \leq n/2$, $|S_k| =  |S_{(n/2)+k}| = 14$, and $d_1(x_a,x_b) = n/2$ for $x_a \in S_k$ and $x_b \in S_{(n/2)+k}$. Hence, by Lemma \ref{lemdist} and Lemma \ref{lem3}, $J_f(P_n \Box H) = J_f(x_0,\ldots,x_{14n-1}) = \sum_{t=0,2,...;t \neq 7n-2}^{14n-4} J_f(x_t,x_{t+1},x_{t+2}) + J_f(x_{14n-2},x_{14n-1}) \geq \sum_{t=0,2,...;t \neq 7n-2}^{14n-4} J_f(x_t,x_{t+1},x_{t+2}) \geq |V(P_n \Box H)|/2-2 = 7n-2$. \qed
\end{proof}

Let $f$ be any radio labeling of $P_n \Box H$, and let $\vec{V}_f= (x_0,x_1,\ldots,x_{14n-1})$ be the ordering of $V(P_n \Box H)$ induced by $f$. Define
$$
\sigma(f) := \sum_{t=0}^{14n-2} \left[ J_f(x_t,x_{t+1})+2\phi(x_t,x_{t+1})-\delta(x_t,x_{t+1})-d_2(x_t,x_{t+1})\right].
$$

\begin{lemma}\label{lem4} Let $f \in \R_o$ be any radio labeling of $P_n \Box H$, and let $\vec{V}_f=(x_0,\ldots,x_{14n-1})$ be the ordering of $V(P_n \Box H)$ induced by $f$. Then
\begin{equation}\label{eqn:sigma1}
\sigma(f) \geq \left\{
\begin{array}{ll}
3-42n, & \mbox{ if $n$ is odd}, \\ [0.2cm]
3-49n, & \mbox{ if $n$ is even}.
\end{array}
\right.
\end{equation}
\end{lemma}
\begin{proof} We consider the following two cases.

\textsf{Case-1:} $n$ is odd. 

Since $f \in \R_o$, we have $\phi(x_t,x_{t+1}) = \delta(x_t,x_{t+1}) = 0$ for all $0 \leq t \leq 14n-2$. Hence, we obtain $\sigma(f) = \sum_{t=0}^{14n-2}(J_f(x_t,x_{t+1})-d_2(x_t,x_{t+1}))$. As $J_f(x_t,x_{t+1}) \geq 0$ and $0 \leq d_2(x_t,x_{t+1}) \leq 3$ for all $0 \leq t \leq 14n-2$, it follows that $\sigma(f) \geq -3(14n-1) = 3-42n$.

\textsf{Case-2:} $n$ is even. 

Since $f \in \R_o$, we have $\phi(x_t,x_{t+1}) = 0$ and $\delta(x_t,x_{t+1}) = 1$ for all $0 \leq t \leq 14n-2$. Hence, we obtain $\sigma(f) = \sum_{t=0}^{14n-2}(J_f(x_t,x_{t+1})-d_2(x_t,x_{t+1})-1)$. Assume that $\vec{V}_f$ is proper. Then, as $f \in \R_o$, we have $J_f(P_n \Box H) = \sum_{t=0}^{14n-1}J_f(x_t,x_{t+1}) \geq 7n-2$ by Lemma \ref{lemjG} and $\sum_{t=0}^{14n-2}d_2(x_t,x_{t+1}) = 3(14n-2)+2$ by Lemma \ref{lem:max}. Hence, we obtain $\sigma(f) \geq 3-49n$. Assume that $\vec{V}_f$ is semi-proper. Then, as $f \in \R_o$, we have $J_f(P_n \Box H) = \sum_{t=0}^{14n-1}J_f(x_t,x_{t+1}) \geq 0$ as $J_f(x_t,x_{t+1}) \geq 0$ for $0 \leq t \leq 14n-2$ and $\sum_{t=0}^{14n-2}d_2(x_t,x_{t+1}) = 3(7n)+2(7n-1)$. Hence, we obtain $\sigma(f) \geq 3-49n$. Assume $\vec{V}_f$ is neither proper nor semi-proper. For any $U_l=(x_l,x_{l+1}x_{l+2})$, denote $\sigma(U_l)=J_f(x_l,x_{l+1},x_{l+2})-d_2(x_l,x_{l+1})-d_2(x_{l+1},x_{l+2})-2$ and $\ell(U_l) = \ell(x_l,x_{l+1},x_{l+2})$. Let $\C = \{U_l : U_l = (x_l,x_{l+1},x_{l+2}), l=0,2,\ldots,14n-4\}$. Then note that $\sigma(f) = \sum_{U_l \in \C}\sigma(U_l)+J_f(x_{14n-2},x_{14n-1})-d_2(x_{14n-2},x_{14n-1})-1$. Let $A = \{U_l \in \C : U_l$ is proper and $\ell(x_l,x_{l+1},x_{l+2}) = n/2\}$, $B = \{U_l \in \C : U_l$ is semi-proper and $\ell(x_l,x_{l+1},x_{l+2}) = n/2\}$ and $C = \{U_l \in \C : U_l$ is neither proper nor semi-proper and  $\ell(x_l,x_{l+1},x_{l+2}) = n/2\}$. If $U_l \in A$, then  $J_f(x_l,x_{l+1},x_{l+2}) \geq 1$ by Lemma \ref{lem3}, $d_2(x_l,x_{l+1})=d_2(x_{l+1},x_{l+2})=3$ and hence $\sigma(U_l) \geq 1-6-2 = -7$. If $U_l \in B$, then $J_f(x_l,x_{l+1},x_{l+2}) \geq 0$ by Lemma \ref{lem3}, $d_2(x_l,x_{l+1})+d_2(x_{l+1},x_{l+2}) = 5$ and hence $J_f(U) \geq 0-5-2=-7$. If $U_l \in C$, then $J_f(x_l,x_{l+1},x_{l+2}) \geq 0$ by Lemma \ref{lem3}, $d_2(x_l,x_{l+1})+d_2(x_{l+1},x_{l+2}) \leq 4$ and hence $\sigma(U_l) \geq 0-4-2 = -6$. If $U_l \in \C \setminus (A \cup B \cup C)$, then note that either $\ell(U_l) = n/2+k$ or $\ell(U_l) = n/2-k$ for some $k \geq 1$. By Lemma \ref{lemdist}, without loss of generality, assume that $\ell(U_l) = n/2+k$ for some $k \geq 1$. By Lemma \ref{lemdist}, there exist $U_{l'}$ such that $\ell(U_{l'})=n/2-k$. Observe that $U_{l'} \in \C \setminus (A \cup B \cup C)$. If $U_l$ is proper, then $\sigma(U_l)+\sigma(U_{l'}) \geq ((2k+1)-6-2)+(0-6-2) \geq -13$. If $U_l$ is semi-proper, then $\sigma(U_l)+\sigma(U_{l'}) \geq ((2k-1)-5-2)+(0-6-2) = -14$. If $U_l$ is weak-proper, then $\sigma(U_l)+\sigma(U_{l'}) \geq ((2k-1)-4-2)+(0-6-2) \geq -13$. If $U_l$ is not proper, semi-proper or weak-proper, then $\sigma(U_l)+\sigma(U_{l'}) \geq (\max\{2k-3,0\}-3-2)+(0-6-2) \geq -13$. Hence $\sigma(f) =  \sum_{U_l \in A}\sigma(U_l) + \sum_{U_l \in B}\sigma(U_l) + \sum_{U_l \in C}\sigma(U_l) + \sum_{U_l \in \C \setminus (A \cup B \cup C)}\sigma(U_l) + J_f(x_{14n-2},x_{14n-1})-d_2(x_{14n-2},x_{14n-1})-1 \geq (|A| + |B| + |C| + |\C \setminus (A \cup B \cup C)|)(-7) - 4 = (7n-1)(-7)-4 = 3-49n$. \hfill \qed
\end{proof}

Let $f$ be any radio labeling of $P_n \Box H$, and let $\vec{V}_f = (x_0,x_1,\ldots,x_{14n-1})$ be the ordering of $V(P_n \Box H)$ induced by $f$. A pair of consecutive vertices $x_t, x_{t+1}$ in the ordering $\vec{V}_f$ is called a {\it bad pair} if $x_t$ and $x_{t+1}$ are on the same side of $P_n \Box H$. Denote by $Z_{f}$ the set of all bad pairs with respect to $\vec{V}_f$, and by $Z_f^c$ the subset of $Z_f$ consisting of those bad pairs $x_t, x_{t+1}$ such that $x_t \in S_{(n+1)/2}$,  or $x_{t+1} \in S_{(n+1)/2}$, or both $x_t, x_{t+1} \in S_{(n+1)/2}$ if $n$ is odd, and $x_t \in S_{n/2} \cup S_{n/2+1}$, or $x_{t+1} \in S_{n/2} \cup S_{n/2+1}$, or both $x_t, x_{t+1} \in S_{n/2} \cup S_{n/2+1}$ if $n$ is even. Recall that $P_n \Box H$ has three sides when $n$ is odd and two sides when $n$ is even.

\begin{observation}\label{rem1} Let $f \in \R_s$ be any radio labeling of $P_n \Box H$. Then  
\begin{equation*}
|Z_{f}| \geq \left\{
    \begin{array}{ll}
    2, & \text{if $n$ is odd}, \\ [0.2cm]
    4, & \text{if $n$ is even},
    \end{array}
    \right.
\end{equation*}
Moreover, when $n$ is even, $|Z_{f}| \equiv 0$ or $2$ (mod $4$); that is, $|Z_f| = 4k$ or $4k+2$ for some $k \geq 1$. 
\end{observation}

\begin{lemma}\label{lem:Js} Let $f \in \R_s$ be a radio labeling of $P_n \Box H$, and let $\vec{V}_f:=(x_0,x_1,\ldots,x_{14n-1})$ be the ordering of $V(P_n \Box H)$ induced by $f$. Then
\begin{equation*}\sigma(f) \geq \left\{
\begin{array}{ll}
5-42n, & \mbox{if $n$ is odd}, \\ [0.2cm]
5-49n, & \mbox{if $n$ is even}.
\end{array}
\right.
\end{equation*}
\end{lemma}
\begin{proof} We consider the following two cases.

\textsf{Case-1:} $n$ is odd. In this case, observe that $\delta(x_t,x_{t+1}) = 0$ for all $0 \leq t \leq 14n-2$. Thus, $\sigma(f) = \sum_{t=0}^{14n-2}(J_f(x_t,x_{t+1})+2\phi(x_t,x_{t+1})-d_2(x_t,x_{t+1}))$. Note that $\phi(x_t,x_{t+1}) \geq 1$ if $x_t,x_{t+1} \in Z_f \setminus Z_f^c$, and $\phi(x_t,x_{t+1}) = 0$ otherwise. If $|Z_f\setminus Z_f^c| \geq 2$, then since $J_f(x_t,x_{t+1}) \geq 0$ and $d_2(x_t,x_{t+1}) \leq 3$ for all $0 \leq t \leq 14n-2$, we obtain $\sigma(f) \geq |Z_f|-3(14n-1) \geq 5-42n$. If $|Z_f \setminus Z_f^c| < 2$, then observe that $|Z_f^c| \geq 2$. If $\vec{V}_f$ is proper, then note that there exist $x_l \in X_1$ or $X_3$ such that $\ell(x_{l-1},x_l,x_{l+1}) \geq \lfloor n/2 \rfloor + 1$. Hence, $J_f(x_{l-1},x_l,x_{l+1}) \geq 2$, and so $\sigma(f) \geq 2-3(14n-1) = 5-42n$. If $\vec{V}_f$ is semi-proper, then note that there exist $x_l \in X_1$ or $X_3$ such that $\ell(x_{l-1},x_l,x_{l+1}) \geq \lfloor n/2 \rfloor + 1$, and hence, $J_f(x_{l-1},x_l,x_{l+1}) \geq 1$. Also, for some $l'$, $d_2(x_{l'},x_{l'+1})+d_2(x_{l'+1},x_{l'+2}) = 5$. Hence, we obtain $\sigma(f) \geq 1-3(14n-3)-5 = 5-42n$. If $\vec{V}_f$ is neither proper nor semi-proper, then for some $l$, $d_2(x_l,x_{l+1})+d_2(x_{l+1},x_{l+2}) \leq 4$. Hence, we obtain $\sigma(f) \geq 0-3(14n-3)-4 = 5-42n$. Hence, in all cases, we obtain $\sigma(f) \geq 5-42n$.   

\textsf{Case-2:} $n$ is even. In this case, observe that $\phi(x_t,x_{t+1}) \geq 1$ if $x_t,x_{t+1} \in Z_f \setminus Z_f^c$, and $\phi(x_t,x_{t+1}) = 0$ if $x_t, x_{t+1} \in V(P_n \Box H) \setminus Z_f$ or $x_t,x_{t+1} \in Z_f^c$. Also,  $\delta(x_t,x_{t+1}) = 1$ if $x_t,x_{t+1} \in V(P_n \Box H) \setminus Z_f$, and $\delta(x_t,x_{t+1}) = 0$ if $x_t,x_{t+1} \in Z_f$. By Lemma \ref{lem3}, we obtain $\sigma(f) \geq (14n-|Z_f|-2)/2+(Z_f \setminus Z_f^c)+(14n-|Z_f|-1)-3(14n-1)$. Since $|Z_f| \geq 4$ by Observation \ref{rem1}, we obtain $\sigma(f) \geq 5-49n$. \qed
\end{proof}

\begin{lemma}\label{lem:fRo} Let $f$ be any radio labeling of $P_n \Box H$, and let $\vec{V}_f := (x_0,x_1,\ldots,x_{14n-1})$ be the ordering of $V(P_n \Box H)$ induced by $f$. If $x_0,x_{14n-1} \in S_{(n+1)/2}$ and $\sigma(f) = 5-42n$, then $f \in \R_o$.     
\end{lemma}
\begin{proof} Assume, to the contrary, that
$x_0,x_{14n-1} \in S_{(n+1)/2}, \sigma(f)=5-42n$, and $f\in \R_s$. Since $f\in \R_s$, it follows that $|Z_f| \geq 2$.

\textsf{Case 1:} $|Z_f \setminus Z_f^{c}| \ge 2$.  
Then there exist indices $k_1,k_2$ such that $\phi(x_{k_1},x_{k_1+1}) \geq 1$, $\phi(x_{k_2},x_{k_2+1}) \geq 1$. Since $J_f(x_t,x_{t+1}) \geq 0$ and $d_2(x_t,x_{t+1}) \leq 3$ for all
$0\leq t \leq 14n-2$, we have
\begin{align*}
\sigma(f)
& = \sum_{t=0}^{14n-2}\left(J_f(x_t,x_{t+1}) + 2\phi(x_t,x_{t+1}) - d_2(x_t,x_{t+1})\right) \\
&\geq 2\left(\phi(x_{k_1},x_{k_1+1}) + \phi(x_{k_2},x_{k_2+1})\right)
   - \sum_{t=0}^{14n-2} d_2(x_t,x_{t+1}) \\
&\ge 2(1+1) - 3(14n-1)
  = 7 - 42n,
\end{align*}
which contradicts $\sigma(f)=5-42n$.

\textsf{Case 2:} $|Z_f \setminus Z_f^{c}| = 1$.  
Since $|Z_f| \geq 2$, we have $|Z_f^{c}| \geq 2$.
Because $|Z_f\setminus Z_f^{c}| = 1$, there exists an index $k$ such that $\phi(x_k,x_{k+1}) \geq 1$. Since $|Z_f^{c}| \geq 2$, there exists an index $a$ such that $\{x_{a-1},x_a,x_{a+1}\} \cap Z_f = \emptyset$, $\ell(x_{a-1},x_a,x_{a+1}) \geq \frac{n+1}{2}$,
and hence $J_f(x_{a-1},x_a,x_{a+1}) \geq 2$.
Thus,
\begin{align*}
\sigma(f)
&\geq J_f(x_{a-1},x_a,x_{a+1}) + 2\phi(x_k,x_{k+1})
   - \sum_{t=0}^{14n-2} d_2(x_t,x_{t+1}) \\
&\geq 2 + 2(1) - 3(14n-1)
  = 7 - 42n,
\end{align*}
which contradicts $\sigma(f) = 5-42n$.

\textsf{Case 3:} $|Z_f \setminus Z_f^{c}| = 0$.  
Then $Z_f^c = Z_f$ and, since $|Z_f| \geq 2$, we have $|Z_f^{c}| \geq 2$.
Because $x_0, x_{14n-1} \in S_{(n+1)/2}$, there exist indices $k_1,k_2$ such that $\ell(x_{k_1-1},x_{k_1},x_{k_1+1}) \geq \frac{n+1}{2}$, $\ell(x_{k_2-1},x_{k_2},x_{k_2+1}) \geq \frac{n+1}{2}$.
Hence, $J_f(x_{k_1-1},x_{k_1},x_{k_1+1}) \geq 2$,
$J_f(x_{k_2-1},x_{k_2},x_{k_2+1}) \geq 2$.

Therefore,
\begin{align*}
\sigma(f)
&\geq J_f(x_{k_1-1},x_{k_1},x_{k_1+1})
  + J_f(x_{k_2-1},x_{k_2},x_{k_2+1})
  - \sum_{t=0}^{14n-2} d_2(x_t,x_{t+1}) \\
&\geq 2 + 2 - 3(14n-1)
  = 7 - 42n,
\end{align*}
which again contradicts $\sigma(f)=5-42n$.

In all three cases, assuming $L(x_0)+L(x_{14n-1})=0$ and $f\in \R_s$ implies $\sigma(f) \geq 7-42n$, contradicting $\sigma(f)=5-42n$.
Therefore, if $L(x_0)+L(x_{14n-1}) = 0$ and $\sigma(f) \geq 5-42n$, then necessarily $f \in \R_o$. \qed
\end{proof}

We are now ready to prove our main Theorem \ref{thm:PHL}. First, We determine the exact radio number for $P_n \Box H$, and then we characterize all optimal radio labelings of $P_n \Box H$.

\begin{Proof} We complete the proof of Theorem \ref{thm:PHL} by proving the following three claims.

\noindent\textbf{Claim 1:} The right-hand side of \eqref{rn:PHL} is a lower bound for $\rn(P_n \Box H)$.

\begin{PRoof}~To prove this claim, it suffices to show that any radio labeling of $P_n \Box H$ has span no less than the right-hand side of \eqref{rn:PHL}. Let $f$ be any radio labeling of $P_n \Box H$, and let $\vec{V}_f:=(x_0,x_1,\ldots,x_{14n-1})$ be the ordering of $V(P_n \Box H)$ induced by $f$. In view of \eqref{eqn:j1}, the function $f$ is given by $f(x_0) = 0$ and $f(x_{i+1})=f(x_i)+d+1-d(x_i,x_{i+1})+J_f(x_i,x_{i+1})$ for $0 \leq i \leq 14n-2$. Summing this equality from $i=0$ to $i=14n-2$, we obtain
\begin{eqnarray}
\Span(f) & = & f(x_{14n-1})-f(x_0) \non \\
& = & \sum_{t=0}^{14n-2}(f(x_{t+1})-f(x_t)) \non \\
& = & \sum_{t=0}^{14n-2}(d+1+J_f(x_t,x_{t+1})-d(x_t,x_{t+1})) \non \\
& = & (14n-1)(d+1)+\sum_{t=0}^{14n-2}(J_f(x_t,x_{t+1})-d(x_t,x_{t+1})) \label{eqn:sum1}
\end{eqnarray}

We consider the following two cases based on the parity of $n$.

\textsf{Case-1:} $n$ is odd.

In this case, $\delta(x_t,x_{t+1}) = 0$ for all $0 \leq t \leq 14n-2$. Hence, by Lemma \ref{lem:phi} (e), we have  $d(x_t,x_{t+1})=L(x_t)+L(x_{t+1})-2\phi(x_t,x_{t+1})+d_2(x_t,x_{t+1})$ for $0 \leq t \leq 14n-2$. Substituting this into \eqref{eqn:sum1}, we obtain
\begin{eqnarray}
\Span(f) & = & (14n-1)(d+1)+\sum_{t=0}^{14n-2}(J_f(x_t,x_{t+1})-L(x_t)-L(x_{t+1})+2\phi(x_t,x_{t+1}) \non \\
& & \hfill -d_2(x_t,x_{t+1})) \non   \\
& = & (14n-1)(d+1)-2\sum_{t=0}^{14n-1}L(x_t)+L(x_0)+L(x_{14n-1})+\sum_{t=0}^{14n-2}(J_f(x_t,x_{t+1}) \non \\ 
& & +2\phi(x_t,x_{t+1})-d_2(x_t,x_{t+1})) \non \\
& = & (14n-1)(d+1)-2L(P_n \Box H)+L(x_0)+L(x_{14n-1})+\sum_{t=0}^{14n-2}(J_f(x_t,x_{t+1}) \non \\ 
& & +2\phi(x_t,x_{t+1})-d_2(x_t,x_{t+1})) \non \\
& = & (14n-1)(d+1)-28L(P_n)+L(x_0)+L(x_{14n-1})+\sigma(f).  \label{eqn:sum2}
\end{eqnarray}

Note that in expression \eqref{eqn:sum2}, the first two terms are constant, while the remaining terms depend on the ordering $\vec{V}_f = (x_0,x_1,\ldots,x_{14n-1})$ induced by $f$. Denote $\pi(f) := L(x_0)+L(x_{14n-1})+ \sigma(f)$. Thus, we aim to determine the optimal value of $\pi(f)$ to establish a lower bound for the radio number of $P_n \Box H$. Since $L(x_0) \geq 0, L(x_{14n-1}) \geq 0$, and $\sigma(f) \geq 3-42n$ by Lemmas \ref{lem4} and \ref{lem:Js}, it follows that $\pi(f) \geq 3-42n$. We now show that there does not exist any radio labeling $f$ whose induced ordering $\vec{V}_f = (x_0,x_1,\ldots,x_{14n-1})$ of $V(P_n \Box H)$, which satisfies $\pi(f) = 3-42n$ for $n \geq 3$ and $\pi(f) =4-42n$ for $n > 3$.

Assume that $\pi(f) = 3-42n$. Since $\pi(f) = L(x_0)+L(x_{14n-1})+\sigma(f)$, and $L(x_0) \geq 0, L(x_{14n-1}) \geq 0$ and $\sigma(f) \geq 3-42n$, we have $L(x_0) = L(x_{14n-1}) = 0$ and $\sigma(f) = 3-42n$. By Lemmas \ref{lem4} and \ref{lem:Js}, the equality $\sigma(f) = 3-42n$ implies that $f \in \R_o$. Hence, $x_t$ and $x_{t+1}$ are on opposite sides for all $0 \leq t \leq 14n-2$. Therefore, by Lemma \ref{lem:phi}, we have $\phi(x_t,x_{t+1}) = 0$ for all $0 \leq t \leq 14n-2$. Moreover, observe that $\sigma(f) = 3-42n$ imply that for every $t$, $J_f(x_t,x_{t+1}) = \phi(x_t,x_{t+1}) = 0$, and $d_2(x_t,x_{t+1}) = 3$. Thus, by Lemmas \ref{lem:phi} and \ref{lem3}, the following holds for the ordering $\vec{V}_f = (x_0,x_1,\ldots,x_{14n-1})$ of $V(P_n \Box H)$: (a) $x_0,x_{14n-1} \in S_{(n+1)/2}$, (b) $x_t$ and $x_{t+1}$ are on opposite sides for all $0 \leq t \leq 14n-2$, (c) $\ell(x_t,x_{t+1},x_{t+2}) \leq (n-1)/2$ for all $0 \leq t \leq 14n-3$, and (d) $d_2(x_t,x_{t+1}) = 3$ for all $0 \leq t \leq 14n-2$. Since $x_0, x_{14n-1} \in S_{(n+1)/2}$, and for any $x_a \in S_1$, $x_b \in S_{(n+1)/2}$, and $x_c \in S_n$, we have $d_1(x_a,x_b)=d_1(x_b,x_c)=(n-1)/2$ and $d_1(x_a,x_c)=n-1$, observe that there exists a vertex $x_t \in S_1$ or $S_n$ such that $d_1(x_{t-1},x_t) > (n-1)/2$ and $d_1(x_t,x_{t+1}) > (n-1)/2$. This contradicts with the condition $\ell(x_{t-1},x_t,x_{t+1}) \leq (n-1)/2$ for all $t$. Hence, we obtain $\pi(f) \geq 4-42n$. For $n = 3$, substituting $d = 5$, $L(P_3) = 2$, and $\pi(f) = -122$ into \eqref{eqn:sum2}, we obtain $\rn(P_3 \Box H) \geq 68$.

For $n > 3$, assume that $\pi(f) = 4-42n$. Since $\pi(f) = L(x_0)+L(x_{14n-1})+\sigma(f)$, and $L(x_0) \geq 0, L(x_{14n-1}) \geq 0$ and $\sigma(f) \geq 3-42n$, there are two possibilities: (1) $L(x_0)+L(x_{14n-1}) = 1$ and $\sigma(f) = 3-42n$; or (2) $L(x_0)+L(x_{14n-1}) = 0$ and $\sigma(f) = 4-42n$. In both cases, by Lemmas \ref{lem4} and \ref{lem:Js}, we have $f \in \R_o$. Hence $x_t$ and $x_{t+1}$ are on opposite sides for all $0 \leq t \leq 14n-2$. Therefore, by Lemma \ref{lem:phi}, we have $\phi(x_t,x_{t+1}) = 0$ for all $0 \leq t \leq 14n-2$. We now consider the two cases separately.

In case (1), $L(x_0)+L(x_{14n-1}) = 1$ and $\sigma(f) = 3-42n$ imply that $L(x_0)+L(x_{14n-1})=1$ and $J_f(x_t,x_{t+1}) = 0, d_2(x_t,x_{t+1}) = 3$ for all $0 \leq t \leq 14n-2$. Since $\Span(f)-f$ is a radio labeling whenever $f$ is a radio labeling, without loss of generality, we may assume that $L(x_0) = 0$ and $L(x_{14n-1}) = 1$. Hence, by Lemmas  \ref{lem:phi}, \ref{lem3} and \ref{lem4}, the following holds for the ordering  $\vec{V}_f = (x_0,x_1,\ldots,x_{14n-1})$ of $V(P_n \Box H)$: (a) $x_0 \in S_{(n+1)/2}$, $x_{14n-1} \not\in S_{(n+1)/2}$, and $x_{14n-1} \in N(x_0)$; (b) $x_t$ and $x_{t+1}$ are on opposite sides for all $0 \leq t \leq 14n-2$; (c) $\ell(x_t,x_{t+1},x_{t+2}) \leq (n-1)/2$ for all $0 \leq t \leq 14n-3$; and (d) $d_2(x_t,x_{t+1}) = 3$ for all $0 \leq t \leq 14n-2$. Note that in this case, all vertices of $S_1,S_{(n+1)/2}$ and $S_n$ must be ordered together. Moreover, the vertices immediate before and after each vertex of $S_{(n+1)/2}$ in the ordering $\vec{V}_f = (x_0,x_1,\ldots,x_{14n-1})$ must be from $S_1$ and $S_n$; otherwise, there exists a vertex $x_t \in S_1$ or $S_n$ such that $d_1(x_{t-1},x_{t}) > (n-1)/2$ and $d_1(x_t,x_{t+1}) > (n-1)/2$, which leads to a contradiction. If this is the case, then note that the remaining vertices form only two sides. By Observation \ref{Obs2}, the ordered vertices on each side belong to the same category, and on opposite side belong to  opposite category. Hence, there exists a vertex $x_t$ such that $d_2(x_t,x_{t+1}) = 2$, which is a contradiction. Therefore, in any case, there does not exist an ordering $x_0,x_1,\ldots,x_{14n-1}$ of $V(P_n \Box H)$ such that $L(x_0)+L(x_{14n-1}) = 1$ and $\sigma(f) = 3-42n$, which is a contradiction. Therefore, our assumption that $L(x_0)+L(x_{14n-1}) = 1$ and $\sigma(f) = 3-42n$ is not true. 

In case (2), since $L(x_0) \geq 0, L(x_{14n-1}) \geq 0$, and $\sigma(f) = \sum_{t=0}^{14n-2}(J_f(x_t,x_{t+1})-d_2(x_t,x_{t+1}))$ with $J_f(x_t,x_{t+1}) \geq 0$ and $0 \leq d_2(x_t,x_{t+1}) \leq 3$, there are again two possibilities: (i) $L(x_0)+L(x_{14n-1}) = 0$, $J_f(x_k, x_{k+1}) = 1$ for a fixed $k$, where $0 \leq k \leq 14n-2$, $J_f(x_t,x_{t+1})=0$ for all $0 \leq t \leq 14n-2, t \neq k$, and $d_2(x_t,x_{t+1}) = 3$ for all $0 \leq t \leq 14n-2$; or (ii) $L(x_0)+L(x_{14n-1})=0$, $J_f(x_t,x_{t+1}) = 0$ and $d_2(x_k,x_{k+1}) = 2$ for a fixed $k$, where $0 \leq k \leq 14n-2$, and $d_2(x_t,x_{t+1}) = 3$ for all $0 \leq t \leq 14n-2, t \neq k$.

In case (i), by Lemmas \ref{lem:phi}, \ref{lem3} and \ref{lem4}, observe that the following hold for the ordering $\vec{V}_f = (x_0,x_1,\ldots,x_{14n-1})$ of $V(P_n \Box H)$: (a) $x_{0}, x_{14n-1} \in S_{(n+1)/2}$, (b) $x_t$ and $x_{t+1}$ are on opposite sides, (c) $\ell(x_t,x_{t+1},x_{t+2}) \leq (n-1)/2$ for all $0 \leq t \leq 14n-3$, and (d) $d_2(x_t,x_{t+1}) = 3$ for all $0 \leq t \leq 14n-2$. In this case, all vertices of $S_1, S_{(n+1)/2}$ and $S_n$ must be ordered smooth such that the vertices immediate before and after of each vertex of $S_{(n+1)/2}$ in the ordering $\vec{V}_f = (x_0,x_1,\ldots,x_{14n-1})$ must be from $S_1$ and $S_n$; otherwise, there exist a vertex $x_t \in S_1$ or $S_n$ such that $d_1(x_{t-1},x_t) > (n-1)/2$ and $d_1(x_t,x_{t+1}) > (n-1)/2$, which contradicts with $\ell(x_{t-1},x_t,x_{t+1}) \leq (n-1)/2$. If this is the case,  then there exists a vertex $x_a \in S_1$ or $S_n$ such that $d_1(x_{a-1},x_a) > (n-1)/2$ and $d_1(x_a,x_{a+1}) > (n-1)/2$. This implies that $J_f(x_{a-1},x_a,x_{a+1}) \geq 2$, which leads to a contradiction.

In case (ii), note that the following hold for the ordering $\vec{V}_f = (x_0,x_1,\ldots,x_{14n-1})$ of $V(P_n \Box H)$: (a) $x_{0}, x_{14n-1}  \in S_{(n+1)/2}$, (b) $x_t$ and $x_{t+1}$ are on opposite sides for all $0 \leq t \leq 14n-2$, (c) $\ell(x_t,x_{t+1},x_{t+2}) \leq (n-1)/2$ for all $0 \leq t \leq 14n-3$, and (d) $d_2(x_k,x_{k+1})=2$ and  $d_2(x_t,x_{t+1}) = 3$ for all $0 \leq t \leq k-1$ and $k+1 \leq t \leq 14n-1$. If $x_k \in S_1$ and $x_{k+1} \in S_{(n+3)/2} \cup \ldots \cup S_{n}$, then remaining vertices forms two sides only and hence there exist a pair of vertices $x_s$ and $x_{s+1}$ such that $d_2(x_s,x_{s+1}) = 2$, which is a contradiction. If $x_k \in V(P_n \Box H) \setminus (S_{1} \cup S_{n})$, then there exists a vertex $x_t \in S_1$ or $S_n$ such that $d_1(x_{t-1},x_{t}) > (n-1)/2$ and $d_1(x_t,x_{t+1}) > (n-1)/2$, which is a contradiction with $\ell(x_{t-1},x_t,x_{t+1}) \leq (n-1)/2$. 

Hence, from cases (i) and (ii), our assumption $L(x_0)+L(x_{14n-1})=0$ and $\sigma(f) = 4-42n$ is not true. 

Thus, from cases (1) and (2), our assumption $\pi(f) = 4-42n$ is not true. 

Therefore, we obtain $\pi(f) \geq 5-42n$. Moreover, note that $d = n+2$ and $L(P_n) = (n^2-1)/4$. Substituting these into  \eqref{eqn:sum2}, we obtain
\begin{eqnarray*}
\rn(P_n \Box H) & \geq & (14n-1)(n+3)-28\left(\frac{n^2-1}{4}\right)+1+1-3(14n-1) \\ 
& = & 7n^2-n+9. 
\end{eqnarray*}

\textsf{Case-2:} $n$ is even.

In this case, $d(x_t,x_{t+1}) = L(x_t)+L(x_{t+1})+\delta(x_t,x_{t+1})-2\phi(x_t,x_{t+1})+d_2(x_t,x_{t+1})$ for $0 \leq t \leq 14n-2$. Substituting this into \eqref{eqn:sum1}, we obtain
\begin{eqnarray}
\Span(f) & = & (14n-1)(d+1)+\sum_{t=0}^{14n-2}(J_f(x_t,x_{t+1})-L(x_t)-L(x_{t+1})-\delta(x_t,x_{t+1}) \non \\ & & +2\phi(x_t,x_{t+1})-d_2(x_t,x_{t+1})) \non \\
& = & (14n-1)(d+1)-2\sum_{t=0}^{14n-1}L(x_t)+L(x_0)+L(x_{14n-1})+\sum_{t=0}^{14n-2}(J_f(x_t,x_{t+1}) \non \\ 
& & +2\phi(x_t,x_{t+1})-\delta(x_t,x_{t+1})-d_2(x_t,x_{t+1})) \non \\
& = & (14n-1)(d+1)-2L(P_n \Box H)+L(x_0)+L(x_{14n-1})+\sum_{t=0}^{14n-2}(J_f(x_t,x_{t+1}) \non \\ 
& & +2\phi(x_t,x_{t+1})-\delta(x_t,x_{t+1})-d_2(x_t,x_{t+1})) \non \\ 
& = & (14n-1)(d+1)-28L(P_n)+L(x_0)+L(x_{14n-1})+\sigma(f). \label{eqn:sum3}
\end{eqnarray}

Note that $L(x_0) \geq 0, L(x_{14n-1}) \geq 0$, and $\sigma(f) \geq 3-49n$ by Lemma \ref{lem4} and \ref{lem:Js}. Moreover, recall  that $d = n+2$ and $L(P_n) = \frac{n(n-2)}{4}$. Substituting these all in \eqref{eqn:sum3}, we obtain
\begin{eqnarray*}
\rn(P_n \Box H) & \geq & (14n-1)(n+3)-28\left(\frac{n(n-2)}{4}\right)+3-49n \\
& = & n(7n+6).
\end{eqnarray*}
\end{PRoof}

\noindent\textbf{Claim 2:} The right-hand side of \eqref{rn:PHL} is an upper bound for $\rn(P_n \Box H)$.

\begin{PROOf}~Note that it is enough to give a radio labeling of $P_n \Box H$ whose span equal to the right-hand side of \eqref{rn:PHL}. For this purpose, we first give an ordering $\vec{V} := (x_0,x_1,\ldots,x_{14n-1})$ of $V(P_n \Box H)$, and then define a radio labeling $f$ on it such that the span of $f$ is equal to the right-hand side of \eqref{rn:PHL}. We consider the following two cases based on the parity of $n$.

\textsf{Case-1:} $n$ is odd.

Note that it is easy to give a radio labeling of $P_3 \Box H$ with $\Span(f) = 68$. Hence, we consider $n > 3$. In this case, let $\alpha = \bigl(\begin{smallmatrix}
  1 & 2 & 3 & 4 & 5 & 6 & 7 & 8 & 9 & 10 & 11 & 12 & 13 & 14 \\
  5 & 2 & 13 & 10 & 7 & 4 & 1 & 12 & 9 & 6 & 3 & 14 & 11 & 8
\end{smallmatrix}\bigr)$, $\beta = \bigl(\begin{smallmatrix}
  1 & 2 & 3 & 4 & 5 & 6 & 7 & 8 & 9 & 10 & 11 & 12 & 13 & 14 \\
  1 & 12 & 9 & 6 & 3 & 14 & 11 & 8 & 5 & 2 & 13 & 10 & 7 & 4
\end{smallmatrix}\bigr)$, $\gamma = \bigl(\begin{smallmatrix}
  1 & 2 & 3 & 4 & 5 & 6 & 7 & 8 & 9 & 10 & 11 & 12 & 13 & 14 \\
  10 & 7 & 4 & 1 & 12 & 9 & 6 & 3 & 14 & 11 & 8 & 5 & 2 & 13
\end{smallmatrix}\bigr)$, $\sigma = \bigl(\begin{smallmatrix}
  1 & 2 & 3 & 4 & 5 & 6 & 7 & 8 & 9 & 10 & 11 & 12 & 13 & 14 \\
  12 & 4 & 10 & 2 & 8 & 7 & 13 & 5 & 11 & 3 & 9 & 1 & 14 & 6
\end{smallmatrix}\bigr)$ and $\tau = \bigl(\begin{smallmatrix}
  1 & 2 & 3 & 4 & 5 & 6 & 7  \\
  1 & 8 & 6 & 13 & 4 & 11 & 2
\end{smallmatrix}$
$\begin{smallmatrix}
8 & 9 & 10 & 11 & 12 & 13 & 14 \\
9 & 7 & 14 & 5 & 12 & 3 & 10
\end{smallmatrix}\bigr)$ be five permutations.
Using these five permutations, we first rename $(u_i,v_j)$, where $1 \leq i \leq n, 1  \leq j \leq 14$, as $(a_r,b_s)$ as follows.
\begin{equation*}
(a_r,b_s) = \left\{
\begin{array}{ll}
(u_i,v_{\alpha(j)}), & \mbox{if $i = 1$}, \\ [0.2cm]
(u_i,v_{\sigma(j)}), & \mbox{if $2 \leq i \leq \lfloor n/2 \rfloor$}, \\ [0.2cm]
(u_i,v_{\beta(j)}), & \mbox{if $i = \lceil n/2 \rceil$}, \\ [0.2cm]
(u_i,v_{\tau(j)}), & \mbox{if $\lceil n/2 \rceil+1 \leq i \leq n-1$}, \\ [0.2cm]
(u_i,v_{\gamma(j)}), & \mbox{if $i = n$}.
\end{array}
\right.
\end{equation*}

We now define the ordering $\vec{V} = (x_0,x_1,\ldots,x_{14n-1})$ as follows: Let $x_t := (a_r,b_s)$, where
\begin{equation*}
t := \left\{
\begin{array}{ll}
3s-1, & \mbox{if $r = 1$ and $1 \leq s \leq 13$}, \\ [0.2cm]
14(\lfloor n/2 \rfloor-1)+41,  & \mbox{if $r = 1$ and $s = 14$}, \\ [0.2cm]
14(\lfloor n/2 \rfloor-r)+2s+39,  & \mbox{if $2 \leq r \leq \lfloor n/2 \rfloor$ and $1 \leq s \leq 7$}, \\ [0.2cm]
14(2\lfloor n/2 \rfloor-r-1)+2s+27,  & \mbox{if $2 \leq r \leq \lfloor n/2 \rfloor$ and $8 \leq s \leq 14$}, \\ [0.2cm]
3s-3, & \mbox{if $r = \lceil n/2 \rceil$}, \\ [0.2cm]
14(n-r-1)+2s+40, & \mbox{if $\lceil n/2 \rceil+1 \leq r \leq n-1$ and $1 \leq s \leq 7$}, \\ [0.2cm]
14(3\lfloor n/2 \rfloor-r-1)+2s+26, & \mbox{if $\lceil n/2 \rceil+1 \leq r \leq n-1$ and $8 \leq s \leq 14$}, \\ [0.2cm]
3s-2, & \mbox{if $r = n$}.
\end{array}
\right.
\end{equation*}

Then note that $L(x_0)=0$, $L(x_{14n-1})=1$, and for $0 \leq t \leq 14n-2$, $x_t$ and $x_{t+1}$ are on opposite sides. Moreover, $d_2(x_t,x_{t+1}) = 3$ for all $0 \leq t \leq 14n-1, t \neq 14(\lfloor n/2 \rfloor-1)+41$ and $d_2(x_t,x_{t+1}) = 2$ when $t = 14(\lfloor n/2 \rfloor-1)+41$. Define $f$ by $f(x_0) = 0$ and $f(x_{i+1}) = f(x_i)+d+1-d(x_i,x_{i+1})$ for $0 \leq i \leq 14n-2$.

\textsf{Claim-i:} $f$ is a radio labeling of $P_n \Box H$ and $\Span(f) = 7n^2-n+9$.

We use Lemma \ref{Prop:2} to prove that $f$ is a radio labeling. Let $(x_i,\ldots,x_j)$, where $0 \leq i < j \leq 14n-1$, be the segment of vertices in the ordering $\vec{V} = (x_0,x_1,\ldots,x_{14n-1})$, such that $j-i+\sum_{t=i}^{j-1}(d-d(x_i,x_j)) < d$. Since $\min\{d(x_t,x_{t+1}),d(x_{t+1},x_{t+2})\} \leq (d+3)/2$ for all $0 \leq t \leq 14n-3$, it follows that $j \leq i+2$. If $j = i+1$, then by the definition of $f$, it is clear that $f$ satisfies \eqref{eqn:rn}. Assume $j = i+2$. Note that $\min\{d(x_t,x_{t+1}),d(x_{t+1},x_{t+2})\} \leq (d+3)/2$ for all $0 \leq t \leq (14n-3)$ in the above-defined ordering. If $0 \leq i \leq 38$, then note that $d(x_i,x_j) \geq (d+1)/2$. Hence,  $f(x_{j})-f(x_i) = 2(d+1)-d-(d+3)/2 = (d+1)/2 \geq d+1-d(x_i,x_j)$. If $i=39$, then $d(x_i,x_j) = (d-1)/2$. Hence $f(x_j)-f(x_i) = 2(d+1)-(d+3)/2-(d-1) = (d+1)/2 \geq d+1-d(x_i,x_j)$. If $i=40$, then $d(x_i,x_j) = (d-1)/2$. Hence $f(x_j)-f(x_i) = 2(d+1)-(d-1)-(d+3)/2 = (d+1)/2 \geq d+1-d(x_i,x_j)$. If $41 \leq i \leq 7(n-3)+38$, then $d(x_i,x_j) \geq 2$, and hence $f(x_j)-f(x_i) \geq 2(d+1)-(d+3)/2-(d+3)/2 = d-1 \geq d+1-d(x_i,x_j)$. If $i = 7(n-3)+39$, then $d(x_i,x_j) = (d-3)/2$, and hence $f(x_j)-f(x_i) = 2(d+1)-(d+3)/2-(d-1) = (d-3)/2 \geq d+1-d(x_i,x_j)$. If $i = 7(n-3)+40$, then $d(x_i,x_j) = (d-1)/2$. Hence $f(x_j)-f(x_i) \geq 2(d+1)-(d-1)-(d+3)/2 = (d+3)/2 = d+1-d(x_i,x_j)$. If $i = 7(n-3)+41$, then $d(x_i,x_j) = (d-1)/2$, and hence $f(x_j)-f(x_i) = 2(d+1)-(d+3)/2-(d+3)/2 = d-1 \geq d+1-d(x_i,x_j)$. If $7(n-3)+42 \leq i \leq 14n-3$, then $d(x_i,x_j) \geq 2$, and hence $f(x_j)-f(x_i) \geq 2(d+1)-(d+3)/2-(d+3)/2 = d-1 \geq d+1-d(x_i,x_j)$. Hence, in each cases, $f$ satisfies \eqref{eqn:rn}. Therefore, $f$ is a radio labeling of $P_n \Box H$. Moreover,
\begin{eqnarray*}
\Span(f) & = & f(x_{14n-1})-f(x_0) \\
& = & (14n-1)(d+1)-\sum_{t=0}^{14n-2}d(x_t,x_{t+1}) \\
& = & (14n-1)(d+1)-\sum_{t=0}^{14n-2}(L(x_t)+L(x_{t+1})+d_2(x_t,x_{t+1})) \\
& = & (14n-1)(d+1)-28(L(P_n))+L(x_0)+L(x_{14n-1})+3(14n-2)+2 \\
\end{eqnarray*}
Substituting $d=n+2$, $L(P_n) = n(n-2)/4$, $L(x_0)=0$ and $L(x_{14n-1})=1$, we obtain $\Span(f) = 7n^2-n+9$, which completes the proof of Claim-i.

\textsf{Case-2:} $n$ is even.

In this case, let $\sigma = \bigl(\begin{smallmatrix}
  1 & 2 & 3 & 4 & 5 & 6 & 7 & 8 & 9 & 10 & 11 & 12 & 13 & 14 \\
  1 & 10 & 5 & 14 & 9 & 4 & 13 & 8 & 3 & 12 & 7 & 2 & 11 & 6
\end{smallmatrix}\bigr)$ and $\tau = \bigl(\begin{smallmatrix}
  1 & 2 & 3 & 4 & 5 & 6 & 7 & 8 & 9 & 10 & 11 & 12 & 13 & 14 \\
  8 & 3 & 12 & 7 & 2 & 11 & 6 & 1 & 10 & 5 & 14 & 9 & 4 & 13
\end{smallmatrix}\bigr)$ be two permutations.

Using these two permutations, we first rename $(u_i,v_j)$, where $1 \leq i \leq n, 1 \leq j \leq 14$, as $(a_r,b_s)$ as follows.
\begin{equation*}
(a_r,b_s) = \left\{
\begin{array}{ll}
(u_i,v_{\sigma(j)}), & \mbox{if $1 \leq i \leq n/2$}, \\ [0.2cm]
(u_i,v_{\tau(j)}), & \mbox{if $n/2 \leq i \leq n$}.
\end{array}
\right.
\end{equation*}

We now define an ordering $\vec{V} = (x_0,x_1,\ldots,x_{14n-1})$ as follows: Let $x_t:=(a_r,b_s)$, where
\begin{equation*}
t := \left\{
\begin{array}{ll}
(n/2-r)28+2(s-1), & \mbox{if $1 \leq r \leq n/2$}, \\ [0.2cm]
(n-r)28+2s-1, & \mbox{if $n/2 \leq i \leq n$}.
\end{array}
\right.
\end{equation*}

Then note that $L(x_0) = L(x_{14n-1}) = 0$, and for $0 \leq t \leq 14n-1$, $x_t$ and $x_{t+1}$ are on opposite sides. Moreover, for $0 \leq t \leq 14n-2$, $d_2(x_t,x_{t+1}) = 3$ when $t$ is even and $d_2(x_t,x_{t+1}) = 2$ when $t$ is odd. Define $f$ by $f(x_0) = 0$, and $f(x_{i+1}) = f(x_i)+d+1-d(x_i,x_{i+1})$ for $0 \leq i \leq 14n-2$.

\textsf{Claim-ii:} $f$ is a radio labeling of $P_n \Box H$ and $\Span(f) = n(7n+6)$.

We use Lemma \ref{Prop:2} to prove that $f$ is a radio labeling of $P_n \Box H$. Let $(x_i,\ldots,x_j)$, where $0 \leq i < j \leq 14n-1$, be the segment of vertices of ordering $\vec{V} = (x_0,x_1,\ldots,x_{14n-1})$, such that $j-i+\sum_{t=i}^{j-1}(d-d(x_t,x_{t+1})) < d$. Since $\min\{d(x_t,x_{t+1}),d(x_{t+1},x_{t+2})\} \leq (d+2)/2$ for all $0 \leq t \leq 14n-3$, it is clear that $j \leq i+2$. If $j=i+1$, then by definition of $f$, it is clear that $f$ satisfies \eqref{eqn:rn}. Assume $j=i+2$. It is clear from above defined ordering that $(d+2)/2 \leq \min\{d(x_i,x_{i+1}),d(x_{i+1},x_{i+2})\} \leq (d+4)/2$ and $\max\{d(x_i,x_{i+1}),d(x_{i+1},x_{i+2})\} = (d+4)/2$ for all $0 \leq i \leq 14n-3$. If $\min\{d(x_i,x_{i+1}),d(x_{i+1},x_{i+2})\} = (d+2)/2$ then note that $d(x_i,x_{i+2}) = 3$. Hence, $f(x_j)-f(x_i) = f(x_{i+2})-f(x_i) = 2(d+1)-d(x_i,x_{i+1})-d(x_{i+1},x_{i+2}) = 2(d+1)-(d+2)/2-(d+4)/2 = d-1 \geq d+1-d(x_i,x_j)$. If $\min\{d(x_i,x_{i+1}),d(x_{i+1},x_{i+2})\} = (d+4)/2$ then note that $d(x_i,x_{i+2}) = 4$. Hence, $f(x_j)-f(x_i) = f(x_{i+2})-f(x_i) = 2(d+1)-d(x_i,x_{i+1})-d(x_{i+1},x_{i+2}) = 2(d+1)-(d+4)/2-(d+4)/2 = d-2 \geq d+1-d(x_i,x_j)$. Hence, in both the cases $f$ satisfies \eqref{eqn:rn}. Thus, $f$ is a radio labeling of $P_n \Box H$. Moreover,
\begin{eqnarray*}
\Span(f) & = & f(x_{14n-1})-f(x_0) \\
& = & (14n-1)(d+1)-\sum_{t=0}^{14n-2}d(x_t,x_{t+1}) \\
& = & (14n-1)(d+1)-\sum_{t=0}^{14n-2}(L(x_t)+L(x_{t+1})+\delta(x_t,x_{t+1})+d_2(x_t,x_{t+1}) \\
& = & (14n-1)(d+1)-\[28(L(P_n))-L(x_0)-L(x_{14n-1})+(14n-1)+35n-2\]
\end{eqnarray*}
Substituting $d=n+2$, $L(P_n) = n(n-2)/4$ and $L(x_0) = L(x_{14n-1}) = 0$, we obtain $\Span(f) = n(7n+6)$.
\end{PROOf}

\noindent\textbf{Claim 3:} The  moreover part of Theorem \ref{thm:PHL} characterizes all optimal radio labelings of $P_n \Box H$.

\begin{PROOF}~Necessity: Let $f : V(P_n \Box H) \rightarrow \mathbb{N} \cup \{0\}$ be an optimal radio labeling of $P_n \Box H$, where $n > 3$, and let $\vec{V}_f = (x_0,x_1,\ldots,x_{14n-1})$ be the ordering of $V(P_n \Box H)$ induced by $f$. We consider the following two cases.

\textsf{Case-1:} $n$ is odd.

It follows from $\Span(f) = 7n^2-n+9$ that $\pi(f) = 5-42n$. Recall that $\pi(f) = L(x_0)+L(x_{14n-1})+\sigma(f)$, with $L(x_0), L(x_{14n-1}) \geq 0$ and $\sigma(f) \geq 3-42n$. Note that $\pi(f) = 5-42n$ holds, then there are following three possibilities: (1) $L(x_0)+L(x_{14n-1})=2$ and $\sigma(f) = 3-42n$, (2) $L(x_0)+L(x_{14n-1})=1$ and $\sigma(f) = 4-42n$, (3) $L(x_0)+L(x_{14n-1})=0$ and $\sigma(f) = 5-42n$. We consider each case separately as follows. 

\textsf{Subcase-1:} $L(x_0)+L(x_{14n-1})=2$ and $\sigma(f) = 3-42n$. 

In this case, $\sigma(f) = 3-42n$ then by Lemmas \ref{lem4} and \ref{lem:phi}, $f$ must be in $\R_o$. Hence $x_t$ and $x_{t+1}$ are on opposite sides of $P_n \Box H$. Thus, $\phi(x_t,x_{t+1}) = 0$ for every $0 \leq t \leq 14n-2$. Hence, $\sigma(f) = \sum_{t=0}^{14n-2}(J_f(x_t,x_{t+1})-d_2(x_t,x_{t+1}))$, where $J_f(x_t,x_{t+1}) \geq 0$ and $0 \leq d_2(x_t,x_{t+1}) \leq 3$ for all $0 \leq t \leq 14n-2$. Now $\sigma(f) = 3-42n$ implies that $J_f(x_t,x_{t+1}) = 0$ and $d_2(x_t,x_{t+1}) = 3$ for all $0 \leq t \leq 14n-2$. Therefore, the ordering  $\vec{V}_f = (x_0,x_1,\ldots,x_{14n-1})$ of $V(P_n \Box H)$ induced by $f$ satisfies the following conditions: (1) $L(x_0)+L(x_{14n-1}) = 2$, (2) $x_t$ and $x_{t+1}$ are on opposite sides of $P_n \Box H$ for all $0 \leq t \leq 14n-2$, (3) $\ell(x_t,x_{t+1},x_{t+2}) = (n-1)/2$ for $t=0,2,\ldots,14n-3$, and (4) $d_2(x_t,x_{t+1}) = 3$ for all $0 \leq t \leq 14n-2$. Now, observe that if the  ordering $\vec{V}_f$ is not smooth, then there exists a vertex $x_t$, either from $S_1$ or $S_n$, such that $d_1(x_{t-1},x_{t}) > (n-1)/2$ and $d_1(x_t,x_{t+1}) > (n-1)/2$, which is a contradiction. If $\vec{V}_f$ is smooth, then by Observation \ref{Obs2}, there exist a vertex $x_t$ such that $d_2(x_t,x_{t+1})=2$, which is also a contradiction. Hence, there does not exist a radio labeling $f$ whose induced ordering $\vec{V}_f = (x_0,x_1,\ldots,x_{14n-1})$ satisfies all the conditions of hypothesis, which contradicts our assumption. Therefore, the case $L(x_0)+L(x_{14n-1})=2$ and $\sigma(f) = 3-42n$ is not possible.

\textsf{Subcase-2:} $L(x_0)+L(x_{14n-1})=1$ and $\sigma(f) = 4-42n$.

In this case, $\sigma(f) = 4-42n$ then by Lemmas \ref{lem4} and \ref{lem:phi}, $f$ must be in $\R_o$. Hence $x_t$ and $x_{t+1}$ are on opposite sides of $P_n \Box H$. Thus, $\phi(x_t,x_{t+1}) = 0$ for every $0 \leq t \leq 14n-2$. Hence, $\sigma(f) = \sum_{t=0}^{14n-2}(J_f(x_t,x_{t+1})-d_2(x_t,x_{t+1}))$, where $J_f(x_t,x_{t+1}) \geq 0$ and $0 \leq d_2(x_t,x_{t+1}) \leq 3$ for all $0 \leq t \leq 14n-2$. Since $\sigma(f) = \sum_{t=0}^{14n-2}(J_f(x_t,x_{t+1})-d_2(x_t,x_{t+1}))$, we again have two possibilities: (i) $L(x_0)+L(x_{14n-1})=1$, $J_f(x_k,x_{k+1})=1$ for a fixed  $k$, $J_f(x_t,x_{t+1})=0$ for all $0 \leq t \leq 14n-2, t \neq k$ and $d_2(x_t,x_{t+1})=3$ for all $0 \leq t \leq 14-2$, and (ii) $L(x_0)+L(x_{14n-1})=1$, $J_f(x_t,x_{t+1})=0$ for all $0 \leq t \leq 14n-2$, $d_2(x_k,x_{k+1})=2$ for a fixed $k$ and $d_2(x_t,x_{t+1}) = 3$ for all $0 \leq t \leq 14n-2, t \neq k$.

We show below that (i) leads to a contradiction.

Assume that (i) is true. Without loss of generality, we may assume that $x_0 \in S_{(n+1)/2}$ and $x_{14n-1} \in N(x_0), x_{14n-1} \notin S_{(n+1)/2}$. Since $d_2(x_t,x_{t+1})=3$ and $J_f(x_t,x_{t+1}) < 2$ for all $0 \leq t \leq 14n-2$, by Lemma \ref{lem3} we have $\ell(x_t,x_{t+1},x_{t+2}) = (n-1)/2$ for $t=0,2,\ldots,14n-3$. Hence, an ordering $\vec{V}_f = (x_0,x_1,\ldots,x_{14n-1})$ of $V(P_n \Box H)$ induced by $f$ satisfies the following conditions: (1) $x_0 \in S_{(n+1)/2}$ and $x_{14n-1} \in N(x_0), x_{14n-1} \notin S_{(n+1)/2}$, (2) $x_t$ and $x_{t+1}$ are on opposite sides of $P_n \Box H$ for all $0 \leq t \leq 14n-2$, (3) $\ell(x_t,x_{t+1},x_{t+2}) = (n-1)/2$ for $t=0,2,\ldots,14n-3$, (4) $d_2(x_t,x_{t+1}) = 3$ for all $0 \leq t \leq 14n-2$. Then, note that the vertices of $S_1, S_{(n+1)/2}$ and $S_n$ are ordered as smooth triples in the ordering $\vec{V}_f$; otherwise there exists a vertex  $x_t$ either from $S_1$ or $S_n$ such that $d_1(x_{t-1},x_{t}) > (n-1)/2$ and $d_1(x_t,x_{t+1}) > (n-1)/2$ which is a contradiction. Observe that $V(P_n \Box H) \setminus (S_1 \cup S_{(n+1)/2} \cup S_n)$ has two sides only and by Observation \ref{Obs2}, the vertices ordered on each side are of same category and hence there exist a vertex $x_s$ to complete the ordering $\vec{V}_f$ such that $d_2(x_s,x_{s+1})=2$, which is a contradiction. Hence there does not exists a radio labeling $f$ satisfying all condition of hypothesis, which contradicts with our assumption.

\textsf{Subcase-3:} $L(x_0)+L(x_{14n-1})=0$ and $\sigma(f) = 5-42n$. 

Since $L(x_0)+L(x_{14n-1}) = 0$, we have $x_0, x_{14n-1} \in S_{(n+1)/2}$. Moreover, since $\sigma(f) = 5-42n$, Lemma \ref{lem:fRo} implies that $f \in \R_o$. Since $f \in \R_o$, $x_t$ and $x_{t+1}$ are on opposite sides of $P_n \Box H$ for all $0 \leq t \leq 14n-2$. Hence, $\phi(x_t,x_{t+1}) = 0$ for all $0 \leq t \leq 14n-2$. Consequently, $\sigma(f) = \sum_{t=0}^{14n-2}(J_f(x_t,x_{t+1})-d_2(x_t,x_{t+1}))$.

Since $\sigma(f) = \sum_{t=0}^{14n-2}(J_f(x_t,x_{t+1})-d_2(x_t,x_{t+1})) = 5-42n$ with $J_f(x_t,x_{t+1}) \geq 0$ and $d_2(x_t,x_{t+1}) \leq 3$ for all $0 \leq t \leq 14n-2$, we have the following three possibilities: (i) $L(x_0)+L(x_{14n-1})=0$, $J_f(x_k,x_{k+1})=2$ for a fixed $k$, $J_f(x_t,x_{t+1}) = 0$ for all $0 \leq t \leq 14n-2, t \neq k$ and $d_2(x_t,x_{t+1})=3$ for all $0 \leq t \leq 14n-2$, (ii) $L(x_0)+L(x_{14n-1})=0$, $J_f(x_k,x_{k+1})=1$ for a fixed $k$, $J_f(x_t,x_{t+1}) = 0$ for all $0 \leq t \leq 14n-2, t \neq k$ and $d_2(x_s,x_{s+1})=2$ for a fixed $s$, $d_2(x_t,x_{t+1})=3$ for all $0 \leq t \leq 14n-2, t \neq s$, and (iii) $L(x_0)+L(x_{14n-1})=0$, $J_f(x_t,x_{t+1})=0$ for all $0 \leq t \leq 14n-2$ and $d_2(x_k,x_{k+1})=d_2(x_s,x_{s+1})=2$ for fixed $k,s$, $d_2(x_t,x_{t+1})=3$ for all $0 \leq t \leq 14n-2, t \neq k,s$. 

We show below that (i) and (ii) leads to a contradiction.. 

Assume that (i) is true. Hence, the ordering $\vec{V}_{f} = (x_0,x_1,\ldots,x_{14n-1})$ of $V(P_n \Box H)$ induced by $f$ satisfies the following: (1) $x_0,x_{14n-1} \in S_{(n+1)/2}$, (2) $x_t$ and $x_{t+1}$ are on opposite sides of $P_n \Box H$ for all $0 \leq t \leq 14n-2$, (3) $\ell(x_k,x_{k+1},x_{k+2}) = (n+1)/2$ and $\ell(x_t,x_{t+1},x_{t+2}) = (n-1)/2$ for all $t = 0,2,\ldots,14n-3, t \neq k$, and (4) $d_2(x_t,x_{t+1})=3$ for all $0 \leq t \leq 14n-2$. Since $d_2(x_t,x_{t+1}) = 3$ for all $0 \leq t \leq 14n-2$, the ordering $\vec{V}_f = (x_0,x_1,\ldots,x_{14n-1})$ of $V(P_n \Box H)$ induced by $f$ is proper. If $\vec{V}_f$ is smooth, then observe that $V(P_n \Box H) \setminus (S_1 \cup S_{(n+1)/2} \cup S_n)$ has two sides only and by Lemma \ref{Obs2}, the vertices ordered on each side are of the same category. Hence, there exists a vertex $x_s$ such that $d_2(x_s,x_{s+1}) = 2$ to complete the ordering $\vec{V}_f$, which is a contradiction. If $\vec{V}_f$ is not smooth, then note that either $L(x_0) = 1$ or $L(x_{14n-1})=1$, which is a contradiction. Hence, there does not exist an ordering $\vec{V}_f$ satisfying all the conditions of hypothesis, which contradicts assumption (i). Therefore, assumption (i) is not true. 

Assume that (ii) is true. Hence, the ordering $\vec{V}_f = (x_0,x_1,\ldots,x_{14n-1})$ of $V(P_n \Box H)$ induced by $f$ satisfies the following: (1) $x_0,x_{14n-1} \in S_{(n+1)/2}$, (2) $x_t$ and $x_{t+1}$ are on opposite sides of $P_n \Box H$ for all $0 \leq t \leq 14n-2$, (3) $\ell(x_k,x_{k+1},x_{k+2}) \leq (n-1)/2$ and $J_f(x_t,x_{t+1}) = 0$ for all $0 \leq t \leq 14n-1, t \neq k, k+1$, and (4) $d_2(x_s,x_{s+1})=2$ and $d_2(x_t,x_{t+1})=3$ for all $0 \leq t \leq 14n-1, t \neq s$. In this case, note that either $L(x_0) = 1$ or $L(x_{14n-1}) = 1$, which is a contradiction. Hence, there does not exist an ordering $\vec{V}_f$ satisfying all conditions of hypothesis, which contradicts assumption (ii). Therefore, assumption (ii) is not true. 

In summary, from above cases, if a function $f$ is an optimal radio labeling of $V(P_n \Box H)$, then the ordering  $\vec{V}_f = (x_0,x_1,\ldots,x_{14n-1})$ induced by $f$ has two possibilities. Hence, by Lemma \ref{lem3} and \ref{lem4}, these can be summarized as follows: (1) $x_t$ and $x_{t+1}$ are on opposite side of $V(P_n \Box H)$ for all $0 \leq t \leq 14n-2$, (2) either $x_0 \in S_{(n+1)/2}$ and $x_{14n-1} \in N(x_0), x_{14n-1} \notin S_{(n+1)/2}$ and $(x_t,x_{t+1},x_{t+2})$ is proper with $\ell(x_t,x_{t+1},x_{t+2}) = \lfloor n/2 \rfloor$ for all $t=0,2,\ldots,(14n-3)/2$; or $x_0,x_{14n-1} \in S_{(n+1)/2}$ and $(x_t,x_{t+1},x_{t+2})$ is proper except two triple which are semi-proper and $\ell(x_t,x_{t+1},x_{t+2}) = \lfloor n/2 \rfloor$. Moreover, by taking $J_t:=J_f(x_t,x_{t+1})$, it is easy to see that $\Span(f) = 7n^2-n+9$ then $f$ must be defined by $f(x_0) = 0$ and $f(x_{t+1}) = f(x_t)+d+1+J_t-L(x_t)-L(x_{t+1})-d_2(x_t,x_{t+1})$, where $J_t$ are as given in \eqref{eqn:Jt}.

\textsf{Case-2:} $n$ is even.

In this case, $\Span(f) = n(7n+6)$. Let $\vec{V}_f = (x_0,x_1,\ldots,x_{14n-1})$ be the ordering of $V(P_n \Box H)$ induced by $f$.  By the proof of Claim 1 and 2, it is clear that $\Span(f) = n(7n+6)$ holds, then $\pi(f) = 3-49n$. Again, if $\pi(f) = 3-49n$, then we obtain $L(x_0)+L(x_{14n-1}) = 0$ and $\sigma(f) = 3-49n$. By the proof of Lemma \ref{lem4} these conditions can be summarized as follows: (1) $x_0 \in S_{n/2}$ and $x_{14n-1} \in S_{n/2+1}$, (2) $x_t$ and $x_{t+1}$ are on opposite side for $0 \leq t \leq 14n-2$, (2) $(x_t,x_{t+1},x_{t+2})$ is proper or semi-proper with $\ell(x_t,x_{t+1},x_{t+2}) = n/2$ or there are equal number of triples $(x_t,x_{t+1},x_{t+2})$ such that $(x_t,x_{t+1},x_{t+2})$ is proper with $\ell(x_{t},x_{t+1},x_{t+2}) = n/2-1$ and $(x_t,x_{t+1},x_{t+2})$ is semi-proper with $\ell(x_t,x_{t+1},x_{t+2}) = n/2+1$. Set $J_t := J_f(x_t,x_{t+1})$ for $0 \leq t \leq 14n-2$. In view of \eqref{eqn:j1}, $f$ is given by $f(x_0) = 0$ and $f(x_{t+1}) = f(x_t)+d+1-d(x_t,x_{t+1})+J_f(x_t,x_{t+1})$ for $0 \leq t \leq 14n-2$. Since $f$ is optimal radio labeling of $P_n \Box H$ and $\Span(f) = n(7n+6)$, it is clear that $J_f(x_t,x_{t+1},x_{t+2}) = 1$ when $(x_t,x_{t+1},x_{t+2})$ is proper with $\ell(x_t,x_{t+1},x_{t+2}) = n/2$. Hence the condition (b) in Theorem \ref{rn:PHL} is satisfied.

Sufficiency: Suppose that an ordering $\vec{V} = (x_0,x_1,\ldots,x_{14n-1})$ of $V(P_n \Box H)$, where $n > 3$, satisfies condition (a) of Theorem \ref{thm:PHL} and $f$ is defined on $x_0,x_1,\ldots,x_{14n-1}$ as in condition (b) of Theorem \ref{thm:PHL}, where $J_0,J_1,\ldots,J_{14n-2}$ are defined in \eqref{eqn:Jt}. We now consider the following two cases.

\textsf{Case-1:} $n$ is odd. It suffices to prove that $f$ is a radio labeling of $P_n \Box H$, and the span of $f$ is the right-hand side of the second line of \eqref{rn:PHL}. Let $x_a$ and $x_b$ be two arbitrary vertices with $0 \leq a < b \leq 14n-1$. If $b-a+\sum_{t=a}^{b-1}(d-d(x_t,x_{t+1})) \geq d$, then, by Lemma \ref{Prop:2}, $f$ is a radio labeling $P_n \Box H$. Assume that $b-a+\sum_{t=a}^{b-1}(d-d(x_t,x_{t+1})) < d$. Since $\ell(x_t,x_{t+1},x_{t+2}) \leq (n-1)/2$ for all $0 \leq t \leq 14n-2$, it follows that $b \leq a+2$. Assume $b=a+1$. Since $x_t$ and $x_{t+1}$ are on opposite sides of $P_n \Box H$, we obtain $f(x_b)-f(x_a) = d+1-L(x_a)-L(x_b)-d_2(x_a,x_b) = d+1-d(x_a,x_b)$. Hence, $f$ is a radio labeling. Assume now that $b=a+2$. If $(x_a,x_{a+1},x_{a+2})$ is proper, then without loss of generality, assume that $d_1(x_a,x_{a+1}) \leq d_1(x_{a+1},x_{a+2})$. Then note that $d(x_a,x_{a+2}) = L(x_{a+2})-L(x_a)+2$. If $\ell(x_a,x_{a+1},x_{a+2}) = (n-1)/2$, then we have $L(x_a)+L(x_{a+1}) = d_1(x_a,x_{a+1}) = (d-3)/2$, and hence we obtain, $L(x_{a+1}) = (d-3)/2-L(x_a)$. Therefore, $f(x_b)-f(x_a) = 2(d+1)-(d(x_a,x_{a+1})+d(x_{a+1},x_{a+2})) = 2(d+1)-(L(x_a)+2(L(x_{a+1}))+L(x_{a+2})+6) = 2(d+1)-(L(x_a)+2((d-3)/2-L(x_a))+L(x_{a+2})+6) = d-1+L(x_a)-L(x_{a+2}) = d+1-(L(x_{a+2})-L(x_a)+2) = d+1-d(x_a,x_{a+2})$. If $\ell(x_a,x_{a+1},x_{a+2}) = (n+1)/2$, then we have $L(x_a)+L(x_{a+1}) = d_1(x_a,x_{a+1}) = (d-1)/2$, and hence $L(x_{a+1}) = (d-1)/2-L(x_a)$. Therefore, $f(x_b)-f(x_a) = 2(d+1)+2-(d(x_a,x_{a+1})+d(x_{a+1},x_{a+2})) = 2(d+1)+2-(L(x_a)+2(L(x_{a+1}))+L(x_{a+2})+6) = 2(d+1)+2-(L(x_a)+2((d-1)/2-L(x_a))+L(x_{a+2})+6) = d-1+L(x_a)-L(x_{a+2}) = d+1-(L(x_{a+2})-L(x_a)+2) = d+1-d(x_a,x_{a+2})$. If $(x_a,x_{a+1},x_{a+2})$ is semi-proper, then $\ell(x_a,x_{a+1},x_{a+2}) = (n-1)/2$. Again, without loss of generality, assume that $d_1(x_a,x_{a+1}) \leq d_1(x_{a+1},x_{a+2})$, then we have $L(x_a)+L(x_{a+1}) = (d-3)/2$, and hence $L(x_{a+1}) = (d-3)/2-L(x_a)$. Observe that $d(x_a,x_{a+2}) = L(x_{a+2})-L(x_a)+3$. Therefore, $f(x_b)-f(x_a) = 2(d+1)-(L(x_a)+2L(x_{a+1})+L(x_{a+2})+5) = 2(d+1)-(L(x_a)+2((d-3)/2)-L(x_a)+L(x_{a+2})+5) = d-L(x_{a+2})-L(x_a) > d+1-(L(x_{a+2})-L(x_a)+3) = d+1-d(x_a,x_{a+2})$. Thus $f$ is a radio labeling of $P_n \Box H$.

Next, we show that in both the cases, the span of $f$ is equal to the right-hand side of \eqref{rn:PHL}. Suppose that conditions (1), (2)-(i), and (3) of (a) hold for the ordering $\vec{V} = (x_0,x_1,\ldots,x_{14n-1})$, and $f$ is defined by (b) in Theorem \ref{thm:PHL}. Then, note that $\sum_{t=0}^{14n-2}J_t = 0$, $\sum_{t=0}^{14n-2}d_2(x_t,x_{t+1}) = 3(14n-2)+2$, and $L(x_0)+L(x_{14n-1})=1$. The span of $f$ is
\begin{eqnarray*}
\Span(f) & = & f(x_{14n-1})-f(x_0) \\
& = & \sum_{t=0}^{14n-2}(f(x_{t+1})-f(x_t)) \\
& = & \sum_{t=0}^{14n-2}(d+1+J_t-d(x_t,x_{t+1})) \\
& = & \sum_{t=0}^{14n-2}(d+1+J_t-L(x_t)-L(x_{t+1})-d_2(x_t,x_{t+1})) \\
& = & (14n-1)(d+1)+\sum_{t=0}^{14n-2}J_t-2\sum_{t=0}^{14n-1}L(x_t)+L(x_0)+L(x_{14n-1})-\sum_{t=0}^{14n-2}d_2(x_t,x_{t+1}) \\
& = & (14n-1)(d+1)-2L(P_n \Box H)+4-42n \\
\end{eqnarray*}
Substituting $d=n+2$ and $L(P_n \Box H) = \frac{7}{2}(n^2-1)$ above, we obtain $\Span(f) = 7n^2-n+9$.

Suppose that conditions (1), (2)-(ii) and (3) of (a) hold for the ordering $\vec{V} = (x_0, x_1, \ldots, x_{14n-1})$, and $f$ is defined by (b) in Theorem \ref{thm:PHL}. Then note that $\sum_{t=0}^{14n-2}J_t = 0$, $\sum_{t=0}^{14n-2}d_2(x_t,x_{t+1}) = 3(14n-3)+4$, and $L(x_0)+L(x_{14n-1})=0$. The span of $f$ is
\begin{eqnarray*}
\Span(f) & = & f(x_{14n-1})-f(x_0) \\
& = & \sum_{t=0}^{14n-2}(f(x_{t+1})-f(x_t)) \\
& = & \sum_{t=0}^{14n-2}(d+1+J_t-d(x_t,x_{t+1})) \\
& = & \sum_{t=0}^{14n-2}(d+1+J_t-L(x_t)-L(x_{t+1})-d_2(x_t,x_{t+1})) \\
& = & (14n-1)(d+1)+\sum_{t=0}^{14n-2}J_t-2\sum_{t=0}^{14n-1}L(x_t)+L(x_0)+L(x_{14n-1})-\sum_{t=0}^{14n-2}d_2(x_t,x_{t+1}) \\
& = & (14n-1)(d+1)-2L(P_n \Box H)+2-42n
\end{eqnarray*}
Substituting $d=n+2$ and $L(P_n \Box H) = \frac{7}{2}(n^2-1)$ in above, we obtain $\Span(f) = 7n^2-n+9$.

\textsf{Case-2:} $n$ is even. It is enough to prove that $f$ is a radio labeling of $P_n \Box H$ and the span of $f$ is the right-hand side of the third line of \eqref{rn:PHL}. Let $x_a$ and $x_b$ be two arbitrary vertices with $0 \leq a < b \leq 14n-1$. If $b-a+\sum_{t=a}^{b-1}(d-d(x_t,x_{t+1})) \geq d$, then, by Lemma \ref{Prop:2}, $f$ is a radio labeling of $P_n \Box H$. Assume that $b-a+\sum_{t=a}^{b-1}(d-d(x_t,x_{t+1})) < d$. Since $\ell(x_t,x_{t+1},x_{t+2}) \leq n/2+1$ for all $0 \leq t \leq 14n-2$, it follows that $b \leq a+2$. Assume $b=a+1$. Since $x_t$ and $x_{t+1}$ are on opposite sides of $P_n \Box H$, we obtain $f(x_b)-f(x_a) = d+1-L(x_a)-L(x_{a+1})-1-d_2(x_a,x_b) = d+1-d(x_a,x_{a+1})$. Hence, $f$ is a radio labeling. Now assume $b=a+2$. Let $(x_a,x_{a+1},x_{a+2})$ be proper. Without loss of generality, assume that $d_1(x_a,x_{a+1}) \leq d_1(x_{a+1},x_{a+2})$. Then note that $d(x_a,x_{a+2}) = L(x_{a+2})-L(x_a)+2$. If $\ell(x_a,x_{a+1},x_{a+2}) = n/2$, then $L(x_a)+L(x_{a+1}) = d_1(x_a,x_{a+1})-1 = (d-2)/2-1 = (d-4)/2$, and hence we have $L(x_{a+1}) = (d-4)/2-L(x_a)$. Therefore,  $f(x_b)-f(x_a) = 2(d+1)+1-(L(x_a)+2L(x_{a+1})+L(x_{a+2})+2+6) = 2d+3-(d+4+L(x_{a+2})-L(x_a)) = d-1-(L(x_{a+2})-L(x_a)) = d+1-(L(x_{a+2})-L(x_a)+2) = d+1-d(x_a,x_{a+2})$. If $\ell(x_a,x_{a+1},x_{a+2}) = n/2-1$, then $L(x_a)+L(x_{a+1}) = d_1(x_a,x_{a+1})-1 = (d-2)/2-1-1 = (d-6)/2$, and hence $L(x_{a+1}) = (d-6)/2-L(x_a)$. Therefore, $f(x_b)-f(x_a) = 2(d+1)-(L(x_a)+2L(x_{a+1})+L(x_{a+2})+2+6) = 2(d+1)-(L(x_a)+2((d-6)/2-L(x_a))+L(x_{a+2})+8) = d-L(x_{a+1})+L(x_a) \geq d+1 - (L(x_{a+1})-L(x_a)+2) = d+1-d(x_a,x_{a+1})$. If $(x_a,x_{a+1},x_{a+2})$ is semi-proper, then $\ell(x_a,x_{a+1},x_{a+2}) = n/2+1$. Again without loss of generality, assume that $d_1(x_a,x_{a+1}) \leq d_1(x_{a+1},x_{a+2})$. Then $L(x_a)+L(x_{a+1}) = d_1(x_a,x_{a+1})-1 = (d-2)/2$, and hence we have $L(x_{a+1}) = (d-2)/2-L(x_a)$. Moreover, note that $d(x_a,x_{a+1}) = L(x_{a+1})-L(x_a)+3$. Hence, $f(x_b)-f(x_a) = 2(d+1)+1-(L(x_a)+2L(x_{a+1})+L(x_{a+2})+2+5) = 2(d+1)+1-(L(x_a)+2((d-2)/2-L(x_a))+L(x_{a+2})+7) = d-2-L(x_{a+2})+L(x_a) = d+1-(L(x_{a+2})-L(x_a)+3) = d+1-d(x_a,x_{a+2})$. Hence, in any case, $f$ is a radio labeling of $P_n \Box H$. Now we show that the span of $f$ is the right-hand side of \eqref{rn:PHL}. Assume that the number of triples $(x_l,x_{l+1},x_{l+2})$ which are semi-proper with $l(x_l,x_{l+1},x_{l+2}) = n/2+1$ are $\beta$. Then, note that $\sum_{t=0}^{14n-2}J_t = 7n-2+\beta$, $\sum_{t=0}^{14n-2}d_2(x_t,x_{t+1}) = 3(14n-\beta-2)-2-2\beta$ and $L(x_0)+L(x_{14n-1})=0$. The span of $f$ is

\begin{eqnarray*}
\Span(f) & = & f(x_{14n-1})-f(x_0) \\
& = & \sum_{t=0}^{14n-2}(f(x_{t+1})-f(x_t)) \\
& = & \sum_{t=0}^{14n-2}(d+1+J_t-d(x_t,x_{t+1})) \\
& = & \sum_{t=0}^{14n-2}(d+1+J_t-L(x_t)-L(x_{t+1})-1-d_2(x_t,x_{t+1})) \\
& = & (14n-1)(d+1)+\sum_{t=0}^{14n-2}J_t-2\sum_{t=0}^{14n-1}L(x_t)+L(x_0)+L(x_{14n-1})-(14n-1) \\
& & -\sum_{t=0}^{14n-2}d_2(x_t,x_{t+1}) \\
& = & (14n-1)d+(7n-2-\beta)-2L(P_n \Box H)-3(14n-\beta-2)-2-2\beta
\end{eqnarray*}
Substituting $d=n+2$ and $L(P_n \Box H) = 7n(n-2)/2$ in above, we obtain $\Span(f) = n(7n+6)$, which completes the proof. 
\end{PROOF}
\end{Proof}

\begin{example} In Table - \ref{tab:odd1}, an ordering $\vec{V} = (x_0,\ldots,x_{69})$ of $P_5 \Box H$ and radio labels of the vertices of $P_5 \Box H$ described in Theorem \ref{thm:PHL} are shown.
\end{example}

\begin{table}[!ht]
    \centering
    \caption{An ordering and it radio labels of vertices of $P_5 \Box H$ with $\rn(P_5 \Box H) = 179$.}\label{tab:odd1}
    \begin{tabular}{|c|cc|cc|cc|cc|cc|cc|}
    \hline
         $(u_i,v_j)\frac{i \rightarrow}{j \downarrow}$ & \multicolumn{2}{|c|}{1} & \multicolumn{2}{|c|}{2} & \multicolumn{2}{|c|}{3} & \multicolumn{2}{|c|}{4} & \multicolumn{2}{|c|}{5}   \\ \hline
        1  & $x_{14} $ & 32  & $x_{65}$ & 167 & $x_{0}$ & 0 & $x_{42}$ & 99 & $x_{28}$ & 66  \\
        2  & $x_{5}$ & 11 & $x_{47}$ & 114 & $x_{33}$ & 77 & $x_{56}$ & 140 & $x_{19}$ & 45  \\
        3  & $x_{38}$ & 88 & $x_{61}$ & 155 & $x_{24}$ & 56 & $x_{52}$ & 129 & $x_{10}$ & 24  \\
        4  & $x_{29}$ & 67 & $x_{43}$ & 102 & $x_{15}$ & 35 & $x_{66}$ & 170 & $x_{1}$ & 3  \\
        5  & $x_{20}$ & 46 & $x_{57}$ & 143 & $x_{6}$ & 14 & $x_{48}$ & 117 & $x_{34}$ & 80  \\
        6  & $x_{11}$ & 25 & $x_{53}$ & 132 & $x_{39}$ & 91 & $x_{62}$ & 158 & $x_{25}$ & 59  \\
        7  & $x_{2}$ & 4 & $x_{67}$ & 173 & $x_{30}$ & 70 & $x_{44}$ & 105 & $x_{16}$ & 38  \\
        8  & $x_{35}$ & 81 & $x_{49}$ & 120 & $x_{21}$ & 49 & $x_{58}$ & 146 & $x_{7}$ & 17  \\
        9  & $x_{26}$ & 60 & $x_{63}$ & 161 & $x_{12}$ & 28 & $x_{54}$ & 135 & $x_{40}$ & 94  \\
        10  & $x_{17}$ & 39 & $x_{45}$ & 108 & $x_{3}$ & 7 & $x_{68}$ & 176 & $x_{31}$ & 73  \\
        11  & $x_{8}$ & 18 & $x_{59}$ & 149 & $x_{36}$ & 84 & $x_{50}$ & 123 & $x_{22}$ & 52  \\
        12  & $x_{55}$ & 137 & $x_{41}$ & 96 & $x_{27}$ & 63 & $x_{64}$ & 164 & $x_{13}$ & 31  \\
        13  & $x_{32}$ & 74 & $x_{69}$ & 179 & $x_{18}$ & 42 & $x_{46}$ & 111 & $x_{4}$ & 10  \\
        14  & $x_{23}$ & 53 & $x_{51}$ & 126 & $x_{9}$ & 21 & $x_{60}$ & 152 & $x_{37}$ & 87  \\ \hline
    \end{tabular}
\end{table}

\begin{example} In Table - \ref{tab:even1}, an ordering $\vec{V} = (x_0,\ldots,x_{83})$ of $P_6 \Box H$ and radio labels of the vertices of $P_6 \Box H$ described in Theorem \ref{thm:PHL} are shown.
\end{example}

\begin{table}[!ht]
    \centering
    \caption{An ordering and it radio labels of vertices of $P_6 \Box H$ with $\rn(P_6 \Box H)= 288$.}\label{tab:even1}
    \begin{tabular}{|c|cc|cc|cc|cc|cc|cc|cc|}
    \hline
         $(u_i,v_j)\frac{i \rightarrow}{j \downarrow}$ & \multicolumn{2}{|c|}{1} & \multicolumn{2}{|c|}{2} & \multicolumn{2}{|c|}{3} & \multicolumn{2}{|c|}{4} & \multicolumn{2}{|c|}{5} & \multicolumn{2}{|c|}{6}  \\ \hline
         1  & $x_{56}$ & 194 & $x_{28}$ & 97 & $x_0$ & 0 & $x_{71}$ & 246 & $x_{43}$ & 149 & $x_{15}$ & 52 \\
         2  & $x_{74}$ & 257 & $x_{46}$ & 160 & $x_{18}$ & 63 & $x_{61}$ & 211 & $x_{33}$ & 114 & $x_{5}$ & 17 \\
         3  & $x_{64}$ & 222 & $x_{36}$ & 125 & $x_{8}$ & 28 & $x_{79}$ & 274 & $x_{51}$ & 177 & $x_{23}$ & 80 \\
         4  & $x_{82}$ & 285 & $x_{54}$ & 188 & $x_{26}$ & 91 & $x_{69}$ & 239 & $x_{41}$ & 142 & $x_{13}$ & 45 \\
         5  & $x_{72}$ & 250 & $x_{44}$ & 153 & $x_{16}$ & 56 & $x_{59}$ & 204 & $x_{31}$ & 107 & $x_{3}$ & 10 \\
         6  & $x_{62}$ & 215 & $x_{34}$ & 118 & $x_{6}$ & 21 & $x_{77}$ & 267 & $x_{49}$ & 170 & $x_{21}$ & 73 \\
         7  & $x_{80}$ & 278 & $x_{52}$ & 181 & $x_{24}$ & 84 & $x_{67}$ & 232 & $x_{39}$ & 135 & $x_{11}$ & 38 \\
         8  & $x_{70}$ & 243 & $x_{42}$ & 146 & $x_{14}$ & 49 & $x_{57}$ & 97 & $x_{29}$ & 100 & $x_{1}$ & 3 \\
         9  & $x_{60}$ & 208 & $x_{32}$ & 111 & $x_{4}$ & 14 & $x_{75}$ & 260 & $x_{47}$ & 163 & $x_{19}$ & 66 \\
         10 & $x_{78}$ & 271 & $x_{50}$ & 174 & $x_{22}$ & 77 & $x_{65}$ & 225 & $x_{37}$ & 128 & $x_{9}$ & 31 \\
         11 & $x_{68}$ & 236 & $x_{40}$ & 139 & $x_{12}$ & 42 & $x_{83}$ & 288 & $x_{55}$ & 191 & $x_{27}$ & 94 \\
         12 & $x_{58}$ & 201 & $x_{30}$ & 104 & $x_{2}$ & 7 & $x_{73}$ & 253 & $x_{45}$ & 156 & $x_{17}$ & 59 \\
         13 & $x_{76}$ & 264 & $x_{48}$ & 167 & $x_{20}$ & 70 & $x_{63}$ & 218 & $x_{35}$ & 121 & $x_{7}$ & 24 \\
         14  & $x_{66}$ & 229 & $x_{38}$ & 132 & $x_{10}$ & 35 & $x_{81}$ & 281 & $x_{53}$ & 184 & $x_{25}$ & 87 \\ \hline
    \end{tabular}
\end{table}

\section*{Concluding Remarks}
In \cite{Das}, Das \emph{et al.} gave a technique to find a lower bound for the radio number of graphs. The result presented in \cite{Das} is an important result about the radio number of arbitrary graphs which can be summarized as follows.

Let $G$ be any simple, finite, connected graph of diameter $d \geq 2$, and let $S \subset V(G)$. Define $N(S)$ as the set of vertices adjacent to vertices of $S$. If $d = 2q$, pick a vertex $v$ of $G$ and set $L_0 = \{v\}$. If $d = 2q+1$, pick a maximal clique $C$ of $G$ and set $L_0 = C$. Now, define $L_1 = N(L_0) \setminus L_0$. Recursively, define $L_{i+1} = N(L_i) \setminus (L_0 \cup L_1 \cup \ldots \cup L_{i})$. These sets are called \emph{layers}, and in particular, $L_i$ is called the $i$-th layer. It is clear that $L_i \neq \phi$ for all $i \in \{0,1,\ldots,q\}$. Let $D_q$ be the subgraph of $G$ induced by the vertices of $L_0 \cup L_1 \cup \ldots \cup L_{q}$. Let $f$ be any optimal radio labeling of $G$, where $\Span(f) = n$. The set $f(V(D_q))$ of colors is called \emph{used colors}, and let the used colors $l(V(D_q))$ be $\{a_1,a_2,\ldots,a_{|D_q|}\}$ such that $a_1 < a_2 < \ldots < a_{|D_q|}$. Using these notations, the authors in \cite{Das} gave the following lower bound for $\rn(G)$.
\begin{equation*}\label{das:lb}
\rn(G) \geq (|D_q|-2q-1)+(\alpha+\beta)+2\left(\sum_{i=0}^{q}|L_i|(q-i)\right), \mbox{ where }a_1 \in L_{\alpha} \mbox{ and } a_{|D_q|} \in L_{\beta}.
\end{equation*}

Since $\alpha+\beta \geq 1$ when $d = 2q$ and $\alpha+\beta \geq 0$ when $d = 2q+1$, we obtain the following lower bound for the radio number of $P_n \Box H$.
\begin{equation}\label{rn:das}
\rn(P_{n} \Box H) \geq \left\{
\begin{array}{ll}
45, & \mbox{ if $n = 3$}, \\ [0.2cm]
7n^2-13n+21, & \mbox{ if $n > 3$ is odd}, \\ [0.2cm]
7n^2-13n+24, & \mbox{ if $n > 3$ is even},
\end{array}
\right.
\end{equation}
Note that the difference between the lower bound for $\rn(P_n \Box H)$ given in \eqref{rn:das} (the right-hand side of \eqref{rn:das}) and $\rn(P_n \Box H)$ given in \eqref{rn:PHL} is $12(n-1)$ when $n > 3$ is odd (for $n = 3$, the difference is $12(n-1)-1 = 23$) and $19(n-1)-5$ when $n > 3$ is even. Since $n$ is an arbitrary positive integer, the difference  between the lower bound for $\rn(P_n \Box H)$ given in \eqref{rn:das} (derived in \cite{Das}) and the exact value of $\rn(P_n \Box H)$ given in \eqref{rn:PHL} (derived in the present work) is arbitrarily large.


\begin{thebibliography}{10}
\bibitem{Bantva1} 
D. Bantva, S. Vaidya, S. Zhou, Radio number of trees, Discrete Appl. Math., 217 (2017), 110--122.

\bibitem{Bantva2}
D. Bantva, D. Liu, Optimal radio labelings of block graphs and line graphs of trees, Theoretical Computer Science, 891 (2021), 90--104. 

\bibitem{Chartrand1} G. Chartrand, D. Erwin, F. Harary, P. Zhang, Radio labelings of graphs, Bull Inst. Combin. Appl., 33 (2001), 77-85.

\bibitem{Chartrand2} G. Chartrand, D. Erwin, P. Zhang, A graph labeling suggested by FM channel restrictions, Bull Inst. Combin. Appl., 43 (2005), 43-57.

\bibitem{Chartrand3} G. Chartrand and P. Zhang, Radio colorings of graphs - A survey, Int. J. Comput. Appl. Math., 2(3) (2007), 237-252.

\bibitem{Das} S. Das, S. Ghosh, S. Nandi and S. Sen, A lower bound technique for radio $k$-coloring, Discrete Math., 340 (2017), 855-861.

\bibitem{Griggs} J. R. Griggs and R. K. Yeh, Labeling graphs with a condition at distance two, SIAM J. Discrete Math., 5 (1992), 586-595.


\bibitem{Kim} B. M. Kim, W. Hwang, B. C. Song, Radio number for the product of a path and a complete graph, J. Comb. Optim., 30 (2015), 139-149.

\bibitem{Liu} D. Liu and X. Zhu, Multi-level distance labelings for paths and cycles, SIAM J. Discrete Math., 19 (2005), 610–621.

\bibitem{West}
D. B. West, Introduction to Graph Theory, Prentice-Hall of India, New Delhi (2001).

\end{thebibliography}
\end{document}